\documentclass[11pt]{amsart}
\usepackage{amssymb}
\usepackage{graphicx}

\input xy
\xyoption{all}
\theoremstyle{plain}
\newtheorem{theorem}{Theorem}[section]
\newtheorem{corollary}[theorem]{Corollary}
\newtheorem{lemma}[theorem]{Lemma}
\newtheorem{axiom}[theorem]{Axiom}

\newtheorem{conjecture}[theorem]{Conjecture}
\newtheorem{proposition}[theorem]{Proposition}
\newtheorem{fact}[theorem]{Fact}
\newtheorem{claim}[theorem]{Claim}

\theoremstyle{definition}
\newtheorem{definition}[theorem]{Definition}
\newtheorem{remark}[theorem]{Remark}

\newtheorem{notation}[theorem]{Notation}

\newtheorem{examples}[theorem]{Examples}
\theoremstyle{remark}

\newcommand{\Ab}{\operatorname{Ab}}
\newcommand{\Hanf}{\operatorname{Hanf}}
\newcommand{\R}{\operatorname{R}}
\newcommand{\id}{\operatorname{id}}
\newcommand{\LS}{\operatorname{LS}}

\newcommand{\dom}{\operatorname{dom}}

\newcommand{\ftp}{\operatorname{tp}}
\newcommand{\lan}{\operatorname{L}}
\newcommand{\gatp}{\operatorname{ga-tp}}
\newcommand{\tp}{\operatorname{ga-tp}}
\newcommand{\gaS}{\operatorname{ga-S}}

\newcommand{\sq}[2]{\sideset{^{#1}}{}{\operatorname{#2}}}

\newcommand{\Mod}{\operatorname{Mod}}

\renewcommand{\phi}{\varphi}

\newcommand{\Union}{\bigcup}
\renewcommand{\P}{\mathcal{P}}
\newcommand{\Pm}{\mathcal{P}^{-}}
\renewcommand{\S}{\mathbf{S}}

\newcommand{\initial}\lessdot
     \newcommand{\infinity}{\infty}
\newcommand{\K}{\operatorname{\mathcal{K}}}

\def\x{\mathbf{x}}

\def\a{\mathbf{a}}
\def\b{\mathbf{b}}
\def\c{\mathbf{c}}

\def\l{\langle}
\def\r{\rangle}

\def\isom{\cong}
\def\submodel{\prec}

\newcommand{\B}{\mathcal B}
\newcommand{\C}{\mathfrak C}

\def\?{?\vadjust
{\vbox to 0pt{\vskip-7pt\hbox to 1.1\hsize{\hfill\huge ?!}}}}

\newbox\noforkbox \newdimen\forklinewidth
\setbox0\hbox{$\textstyle\smile$}
\setbox1\hbox to \wd0{\hfil\vrule width \forklinewidth depth-2pt
     height 10pt \hfil}
\setbox\noforkbox\hbox{\lower 2pt\box1\lower 2pt\box0\relax}
\def\unionstick{\mathop{\copy\noforkbox}\limits}

\def\nonfork{\unionstick}

\setbox0\hbox{$\textstyle\smile$}
\setbox1\hbox to \wd0{\hfil{\sl /\/}\hfil}
\setbox2\hbox to \wd0{\hfil\vrule height 10pt depth -2pt width
                   \forklinewidth\hfil}
\newbox\doesforkbox
\setbox\doesforkbox\hbox{\lower 2pt\box1 \lower
2pt\box2\lower2pt\box0\relax}

\begin{document}

\title{Excellent Abstract Elementary Classes are tame}

\date{September 8, 2005}

\author{Rami Grossberg}
\email[Rami Grossberg]{rami@cmu.edu}
\urladdr{http://www.math.cmu.edu/$\sim$rami}
\address{Department of Mathematics\\
Carnegie Mellon University\\
Pittsburgh PA 15213}

\author{Alexei S. Kolesnikov}
\email{alexeik@umich.edu}
\urladdr{http://www.math.lsa.umich.edu/$\sim$alexeik}
\address{
Department of Mathematics\\
University of Michigan\\
Ann Arbor, MI 48109}

\thanks{The research is part of the author's work towards his
Ph.D. degree under direction of Prof. Rami Grossberg. I am deeply
grateful to him for his guidance and support.}


\begin{abstract}

The assumption that an  AEC is tame is a powerful assumption permitting
development of stability theory for AECs with the amalgamation property.
Lately several upward categoricity theorems were discovered where
tameness replaces strong set-theoretic assumptions.

We present in this article two sufficient conditions for tameness, both
in form of strong amalgamation properties that occur in nature.  One of
them was used recently to prove that several Hrushovski classes are tame.

This is done by introducing the property of weak
$(\mu,n)$-uniqueness which makes sense for all AECs (unlike Shelah's
original property) and derive it from the assumption that
weak $(\LS(\K),n)$-uniqueness, $(\LS(\K),n)$-symmetry and
$(\LS(\K),n)$-existence properties hold for all $n<\omega$.  The
conjunction of  these three properties we call
\emph{excellence}, unlike \cite{Sh 87b} we do not require the very
strong $(\LS(\K),n)$-uniqueness, nor we assume that the members of
$\K$ are atomic models of a countable first order theory.  We also work in
a more general context than  Shelah's good frames.

\end{abstract}

\maketitle

%%%%%%%%%%%%% body %%%%%%%%%%%%%%%%%%%%%%%%%%%%

\section*{Introduction}\label{s:introduction}

In 1977 Shelah influenced by earlier work of J\'{o}nsson (\cite{Jo1} and
\cite{Jo2}) in
\cite{Sh 88} introduced a  semantic generalization of Keisler's \cite{Ke}
treatment of
$L_{\omega_1,\omega}(\mathbf Q)$.  It is the notion of
\emph{Abstract Elementary Class}:
\begin{definition}\label{def AEC}
Let $\K$ be a class of structures
all in the same similarity type $\lan(\K)$, and
let $\prec_{\K}$ be a partial order on $\K$.  The ordered pair
$\langle \mathcal{K},\prec_{\mathcal{K}}\rangle$ is an
\emph{abstract elementary class,  AEC for short}
     iff
\begin{enumerate}

\item [A0] (Closure under isomorphism)
\begin{enumerate}

\item
For every $M\in \mathcal{K}$
     and every $\lan(\K)$-structure $N$ if $M\cong N$ then $N\in
\mathcal{K}$.
\item Let $N_1,N_2\in\K$ and  $M_1,M_2\in \K$ such that
there exist $f_l:N_l\cong M_l$ (for $l=1,2$) satisfying
$f_1\subseteq f_2$ then  $N_1\prec_{\K}N_2$
implies that  $M_1\prec_{\K}M_2$.

\end{enumerate}

\item [A1]
For all $M,N\in \mathcal{K}$ if $M\prec_{\mathcal{K}} N$ then $M\subseteq
N$.
     \end{enumerate}

\begin{enumerate}
\item [A2]
Let $M,N,M^*$ be $\lan(\mathcal{K})$-structures.
If $M\subseteq N$, $M\prec_{\mathcal{K}} M^*$ and $N\prec_{\mathcal{K}}
M^*$ then
     $M\prec_{\mathcal{K}} N$.

\item [A3]
(Downward L\"owenheim-Skolem)
There exists a cardinal\\
$\LS(\mathcal{K})\geq \aleph_0+|\lan(\mathcal{K})|$
%\footnote{The purpose of this requirement is to avoid ``pathologies'',
%as we are not ready to consider situations where the cardinality of the
%set of
% terms in the language is greater than that of the structure.}
such that
for every \\
$M\in \mathcal{K}$
     and for every $A\subseteq |M|$ there exists $N\in \mathcal{K} $ such
that
$N\prec_{\mathcal{K}} M, \; |N|\supseteq A$ and
$\|N\|\leq |A|+\LS(\mathcal{K})$.

\item [A4]
(Tarski-Vaught Chain)
\begin{enumerate}

\item
For every regular cardinal  $\mu$ and every \\
$N\in \mathcal{K}$ if
$\{M_i\prec_{\mathcal{K}} N\;:\;i<\mu\}\subseteq \mathcal{K}$
is $\prec_{\K}$-increasing (i.e. $i<j\Longrightarrow M_i\prec_{\mathcal{K}}
M_j$) then
$\Union_{i<\mu}M_i\in \mathcal{K}$ and
$\Union_{i<\mu}M_i\prec_{\mathcal{K}} N
$.

\item
For every regular $\mu$,
if  $\{M_i\;:\;i<\mu\}\subseteq \mathcal{K}$ is
$\prec_{\K}$-increasing then $\Union_{i<\mu}M_i\in \mathcal{K}$ and
$M_0\prec_{\mathcal{K}} \Union_{i<\mu}M_i$.
\end{enumerate}
\end{enumerate}

For $M$ and $N\in \K$ a monomorphism $f:M\rightarrow N$ is called a
\emph{$\K$-embedding} iff $f[M]\prec_{\K}N$.  Thus, $M\prec_{\K}N$ is
equivalent to ``$\id_M$ is a $\K$-embedding from $M$ into $N$''.
     \end{definition}

Many of the fundamental facts on AECs are due to Saharon Shelah and were
introduced in
\cite{Sh 88},
\cite{Sh 394}
and
\cite{Sh 576}.  For a  survey of some of the basics see  \cite{Gr1} or
\cite{Gr3}.

In the  late seventies Shelah established the program of developing
\emph{Classification Theory for Abstract Elementary Classes}, namely that
there exists a vastly more general theory than the one presented in
\cite{Sh c} that can be developed without any reference to the compactness
theorem (that fails already in small fragments of $L_{\omega_1,\omega}$).
As such a theory undoubtedly will require new concepts and techniques
Shelah proposed the following as a test problem:

\begin{conjecture}[Shelah's conjecture]  Let $\psi\in L_{\omega_1,\omega}$
be a sentence
in a countable language.
If $\psi$ is $\lambda$-categorical in some $\lambda > \beth_{\omega_1}$
then $\psi$ is $\mu$-categorical
for every $\mu\geq\beth_{\omega_1}$.

\end{conjecture}

Several authors  wrote many papers trying to
approximate this conjecture (Shelah alone produced  more than
1,000 pages), the conjecture at present seems to be not accessible.

In 1990 Shelah proposed a generalization for AECs:

\begin{conjecture}[see \cite{Sh c}]  Let $\K$ be an AEC.  If $\K$ is
categorical in some
$\lambda > \Hanf(\K)$ then $\K$ is $\mu$-categorical for every
$\mu\geq\Hanf(\K)$.

\end{conjecture}

\begin{notation}
Let $\mu$ be a cardinal number and $\K$ a class of models.  By $\K_\mu$ we
denote
the subclass $\{M\in \K\;:\; \|M\|=\mu\}$.
\end{notation}

Two classical concepts that introduced in the fifties and studied
extensively by Fraisee, Robinson and Jonsson play also an important role
in AECs:

\begin{definition}
Let $\l\K,\prec_{\K}\r$ be an AEC and suppose $\mu \geq \LS(\K)$.  We say
that \emph{$\K$ has
the $\mu$-amalgamation property} iff for all $M_\ell\in \K_\mu$ (for
$\ell=0,1,2$) such
that $M_0\prec_{\K}M_\ell$ (for $\ell =1,2$) there exists
$N^{*}\in \K_\mu$ and $f_\ell:M_\ell\rightarrow N^*$ (for $\ell=1,2$) such
that
$f_1\restriction M_0=f_2\restriction M_0$, i.e.
the following diagram commutes:
\[
\quad \xymatrix{M_1 \ar[r]^{f_1} & N^{*} \\
M_0 \ar[u]^{\id} \ar[r]_{\id} & M_2 \ar[u]_{f_2}
}
\]

The model $N^*$ is called an \emph{amalgam} of $M_1$ and $M_2$ over $M_0$.

\emph{$\K$ has the $\mu$-joint mapping property}  iff for any $M_\ell\in
\K_\mu$
for $\ell=1,2$ there are $N^*\in \K_\mu$ and $\K$-embeddings
$f_\ell:M_\ell\rightarrow N^*$.

We say that $\K$ has the \emph{amalgamation property}  iff it has the
$\mu$-amalgamation property for all $\mu\geq \LS(\K)$.
\end{definition}

Using the axioms of AECs one can prove the following:

\begin{fact}  If $\K$ has the $\mu$-AP for all $\mu\geq\LS(\K)$ then for
any triple $M_\ell\in \K_{\geq\LS(\K)}$.  If $M_0\prec_{\K}M_1,M_2$ then
there exists an amalgam of $M_1$ and $M_2$ over $M_0$.
\end{fact}

Using Axiom A0 from the definition of AEC it follows that both a
stronger-looking and a
weaker-looking amalgamation properties are equivalent to what we call above
the amalgamation
property:

\begin{fact} \label{AP lemma}Let $\K$ be an AEC.  The following are
equivalent
\begin{enumerate}
\item $\K$ has the $\mu$-amalgamation property,

\item
for all $M_\ell\in \K_\mu$ (for $\ell=0,1,2$) such
that $M_0\prec_{\K}M_\ell$ (for $\ell =1,2$) there exists
$N^{*}\in \K_\mu$ such that $N^{*}\succ_{\K}N_2$ and there is
$f:M_1\rightarrow
N$ satisfying $f\restriction M_0=\id_{M_0}$, i.e.
the following diagram commutes:
\[
\quad \xymatrix{M_1 \ar[r]^{f} & N^{*} \\
M_0 \ar[u]^{\id} \ar[r]_{\id} & M_2 \ar[u]_{\id}
}
\]

\item
for all $M_\ell\in \K_\mu$ (for $\ell=0,1,2$) such
that $g_\ell:M_0\rightarrow M_\ell$ (for $\ell =1,2$) are
$\K$-embeddings
there are
$N^{*}\in \K_\mu$  and there is
$f_\ell:M_\ell\rightarrow
N^{*}$ satisfying $f_1\circ g_1\restriction M_0=f_2\circ g_2\restriction
M_0$ i.e. the next diagram commutes:
\[
\quad \xymatrix{M_1 \ar[r]^{f_1} & N^{*} \\
M_0 \ar[u]^{g_1} \ar[r]_{g_2} & M_2 \ar[u]_{f_2}
}
\]

\end{enumerate}
\end{fact}

An important tool in proving the above lemma is the following basic
property of
AECs.

\begin{fact}\label{pull-back}
Suppose $f:M\to N$ is a $\K$-embedding. There are a model $\bar M \succ M$
and
$\bar f:\bar M \cong N$ extending $f$.
\end{fact}

     Robinson's consistency property implies that if $T$ is a
complete first-order  theory then $\Mod(T)$ has both the amalgamation
     and the joint mapping properties.  As there are natural examples of AECs
where  the $\mu$-AP is a property fails (see \cite{GrSh}) we must deal with
AP as a property.

\emph{Acknowledgment: }  John Baldwin provided us with detailed very helpful comments, remarks and questions  on the 8/30/2005 version that improved very much the presentation of this paper.

\section*{Galois types, amalgamation and tameness}

In the theory of AECs the notion of complete first-order type is
replaced by that of a \emph{Galois type}:

\begin {definition}\label{E}\index{$E$, binary relation}
Let $\beta>0$ be an ordinal.
For triples $(\bar a_\ell, M, N_\ell)$ where $\bar a_\ell\in
\sq{\beta}N_\ell$ and
$M\prec_{\K}N_\ell\in\K$ for
$\ell=1,2$,
we define a binary relation $E$ as follows:
$(\bar a_2, M, N_2)E(\bar a_1, M, N_1)$ iff
     and there exists $N\in\K$ and $\K$-mappings
$f_1, f_2$ such that
$f_\ell:N_l\rightarrow N$ and $f_\ell\restriction M=\id_M$ for $\ell=1,2$
and
$f_2(\bar a_2)=f_1(\bar a_1)$:

\[
\xymatrix{\ar @{} [dr] N_1
\ar[r]_{f_1}  &N \\
M \ar[u]^{\id} \ar[r]_{\id}
& N_2 \ar[u]_{f_2}
}
\]

\end{definition}

\begin{remark}  When $\K$ has the amalgamation property then
$E$ is an equivalence relation on the class of triples of the form
$(\bar a, M, N)$.  If $\K$ fails to have the amalgamation property,
$E$ may fail to be transitive, but the transitive closure of $E$
could be used instead.
\end{remark}

\begin{remark}  Using Ax0 one can show that in the previous definition
we may assume that $f_2=\id_{N_2}$, i.e. that $N\succ_{\K}N_2$ and the
condition is that $f_1(\bar a_1)=\bar a_2$.
\end{remark}

\begin{definition}\label{defn of ga types}\index{Galois-type}\index{type,
Galois} Let $\beta$ be a positive ordinal.
\begin{enumerate}

\item For $M,N\in\K$ and $\bar a\in \sq{\beta}N.$
The \emph{Galois type of $\bar a$ in $N$ over $M$}, written
$\tp(\bar a/M,N)$, is defined to be $(\bar a,M,N)/E$.

%\item We abbreviate $\tp(\bar a/M,N)$ by $\tp(\bar a/M)$.\index{$\tp(\bar
%a/M)$}

\item For $M\in\K$,
$$\gaS^\beta(M):=\{\tp(\bar a/M,N)\mid M\prec
N\in\K_{\|M\|},
\bar a\in \sq{\beta}N\}.$$\index{$\gaS^\beta(M)$}
      We write $\gaS(M)$ for
$\gaS^1(M)$.\index{$\gaS^1(M)$}\index{$\gaS^\beta(M)$}

\item Let $p:=\tp(\bar a/M',N)$ for $M\prec_{\K}M'$ we denote by
$p\restriction M$ the type
$\tp(\bar a/M,N)$.  The \emph{domain of $p$} is denoted by $\dom p$ and
it is by definition $M'$.
\index{Galois-type!restriction}\index{Galois-type!domain}

\item  Let $p=\tp(\bar a/M,N)$, suppose that $M\prec_{\K}N'\prec_{\K}N$
and let $\bar b\in \sq{\beta}N'$ we say that \emph{$\bar b$ realizes $p$}
iff
$\tp(\bar b/M,N')=p\restriction M$.\index{Galois-type!realized}

\item For types $p$ and $q$, we write $p\leq q$ if $\dom(p)\subseteq
\dom(q)$ and there exists $\bar a$ realizing $p$ in some $N$ extending
$\dom(p)$ such that $(\bar a,\dom(p), N) = q\restriction
\dom(p)$.\index{Galois-type!extension}
\end{enumerate}
\end{definition}

An important notion in this paper is that of an \emph{amalgamation base}.  A
model is an amalgamation base iff every pair of models extending it of the
same
cardinality can be amalgamated over it.     Sometimes we will
be interested to consider amalgamation bases some special sets which are
not models.
Please note that the assumption that every subset of a model (from $\K$)
is an amalgamation basis
a very strong assumption.  Making this assumption brings us to the very
special context of AECs called \emph{homogeneous model theory},  see
\cite{GrLe} for an introduction.   Since there are many interesting
examples of
AECs with amalgamation over models (but not over all sets) like in Zilber's
theory of pseudo exponentiation we do not make the assumption that all sets
are amalgamation bases.
Our
interest is limited for very special sets that originate from certain
systems
of models we describe now.

\begin{definition}
Let $I$ be a subset of $\mathcal {P}(n)$ for some $n<\omega$ that is
downward closed
(i.e.
$t\in I$ and $s\subseteq t$ implies $s\in I$).

For an
$\mathbf S=\l M_s\in\K \mid s\in I\r$  is an
\emph{$I$-system}\index{$I$-system} iff for all $s,t\in I$
\begin{enumerate}

\item
$s\subseteq  t\implies M_s\prec_{\K} M_t$ and

\item
$M_{s\cap t}=M_s\cap M_t$
\end{enumerate}
$\mathbf S$ is a \emph{$(\lambda,I)$-system}
\index{$(\lambda,I)$-system}iff in addition all the
models are  of cardinality $\lambda$.

Denote by
\[
A^{\mathbf S}_t:=\bigcup_{s\subsetneq t}M_s  \text{ and }
A^{\mathbf S}_I:=\bigcup_{s\in I}M_s
\]

\end{definition}

\begin{definition}\label{def: AB}  Suppose $\mathbf S=\l M_s\in\K_\mu \mid
s\in I\r$  is an
\emph{$I$-system} for some $I\subseteq \P(n)$
We say that \emph{a set $A^{\mathbf S}_I$ is a $\mu$-amalgamation base} iff
for all $M_\ell\in \K_\mu$ (for $\ell=1,2$) such
that $M_s\prec_{\K}M_\ell$ (for all $s\in I$  and $\ell =1,2$) there exists
$N^{*}\in \K_\mu$ such that $N^{*}\succ_{\K} M_2$ and there is a
$\K$-embedding
$f:M_1\rightarrow N^*$ satisfying $f\restriction {A^{\mathbf
S}_I}=\id_{A^{\mathbf S}_I}$, i.e. the following diagram commutes:
\[
\quad \xymatrix{M_1 \ar[r]^{f} & N^{*} \\
{A^{\mathbf S}_I} \ar[u]^{\id_{A^{\mathbf S}_I}} \ar[r]_{\id_{A^{\mathbf
S}_I}}
& M_2
\ar[u]_{\id_{M_2}} }
\]
\end{definition}

\begin{notation}
Denote by $\Ab_\mu(\K)$ the class $$\{{A^{\mathbf S}_I}\mid {A^{\mathbf
S}_I}
\text{ is a
$\mu$-amalgamation base for some }I\text{-system from  }\K_\mu\}.$$
%If $\mu =|A|$, we simply say that $A$ is an amalgamation base.
\end{notation}

Thus $\K$ has the $\lambda$-amalgamation property iff $ \K_\lambda
\subseteq \Ab_\lambda(\K)$.
Under the assumption that $\K_\mu$ has the AP  the notion of a Galois-type
can be extended to include also
$\gatp(\bar a/A,M)$ for $A\in \Ab_\mu (\K)$.

\begin{definition}
Let $\K$ be an AEC with the amalgamation property
and let $\chi\geq \LS(\K)$.  The class $\K$ is called
\emph{$\chi$-tame} iff \index{tame}
\[
p\neq q\implies \exists N\prec_{\K}M \text{ of cardinality }\leq\chi
\text{ such that }p\restriction N\neq q\restriction N
\]
for any $M\in \K_{>\chi}$ and
every $p,q\in \gaS(M)$

$\K$ is \emph{tame} iff it is $\chi$-tame for some $\chi<\Hanf(\K)$

Suppose $\mu>\chi$.  The class is \emph{$(\chi,\mu)$-tame}
\index{$(\chi,\mu)$-tame}
iff
\[
p\neq q\implies \exists N\prec_{\K}M \text{ of cardinality }\leq\chi
\text{ such that }p\restriction N\neq q\restriction N
\]
for any $M\in \K_{\mu}$ and
every $p,q\in \gaS(M)$
\end{definition}

In \cite{GrV1} Grossberg and VanDieren introduced the notion of
\emph{tameness} as a candidate for a further ``reasonable''
assumption an AEC that permits development of stability-like
theory.   It turns out that essentially the same property was
introduced earlier by Shelah implicitly in the proof of his main
theorem in \cite{Sh 394}.

One of the better approximations to Shelah's categoricity conjecture for
AECs can be derived from a theorem due to Makkai and Shelah
(\cite{Sh285}):

\begin{theorem}[Makkai and Shelah 1990]
     Let $\K$ be an AEC,  {  $\kappa$ a strongly compact cardinal}
such that $\LS(\mathcal{K})<\kappa$.  Let $\mu_0:=\beth_{(2^\kappa)^+}$.
If
$\K$ is categorical in some
     $\lambda^{+}>\mu_0$ then $\K$ is categorical in every $\mu\geq\mu_0$.
\end{theorem}

 Proposition 1.13 of \cite{Sh285} asserts  (using the assumption that $\kappa$ is strongly compact)
that any AEC $\K$ as above has the AP (for models of cardinality
$\ge\kappa$). Since Galois types in this context are sets of $L_{\kappa,\kappa}$ formulas the class is trivially $\kappa$-tame.

     In
\cite{GrV2} Grossberg and VanDieren proved (in ZFC)  a case of  Shelah's
categoricity conjecture for  tame AECs with  the amalgamation property
which implies the above theorem of Makkai and Shelah. Thus the
tameness assumption enables upward categoricity argument (instead of the
large cardinal assumption).  This is also an extension (upward) of Shelah's
main theorem from \cite{Sh 394}.

\begin{theorem}[Grossberg and VanDieren 2003]
    Let $\K$ be an AEC, $\kappa:=\beth_{(2^{\LS(\K)})^+}$.  Denote by
$\mu_0:=\beth_{(2^\kappa)^+}$.  Suppose that $\K_{>\kappa}$ has the
amalgamation property and is tame.
If
$\K$ is categorical in some
     $\lambda^{+}>\mu_0$ then $\K$ is categorical in every $\mu\geq\mu_0$.
\end{theorem}

Later  Lessmann obtained finer upward categoricity
results by using much stronger assumptions to tameness ($\aleph_0$-tameness
and $\LS(\K)=\aleph_0$) and existence of arbitrary large models.

In \cite{Sh 394} Shelah proved that for an AEC with the amalgamation
property.  If
$\K$ is $\lambda$-categorical for some $\lambda>\beth_{(2^{\Hanf(\K)})^+}$
then it is
$(\Hanf(\K),\mu)$-tame for all $\Hanf(\K)<\mu<\lambda$.

Throughout this paper we will be using Shelah's presentation theorem for
AECs which states that
every AEC can be viewed as a PC-class (see \cite{Sh 88} or \cite {Gr3}). We
state it in a form that is
more convenient for our purposes.

\begin{lemma}\label{Skolem}
Let $\K$ be an AEC, let $\mu=\LS(\K)$. Let $\chi_0$ e a large regular
cardinal.
There are $\mu$ functions $\{f_i \mid i<\mu\}$ such that whenever $M\in
\K$, $M\subset H(\chi_0)$,
and $\B \submodel \l H(\chi_0),\in, \K, M, \{f_i\mid i<\mu\} \r$, $\|\B\|\ge
\mu $, for $N=M^{\B}$
we have $N\in \K$ and $N\submodel_{\K} M$.
\end{lemma}

This is simply saying that Skolem functions can be defined in an
appropriate set-theoretic universe and whenever a subset $N$ of a model
$M\in \K$ is closed under those functions, $N$ is a $\K$-model.

%end introduction

%%%%%%%%%%%%%%%%%%%%%%%%%%%%%%%%%%%%%%%%%%%%%%%%%%%%%%%%%%%%%%%%
\section{The basic framework and concepts}
\label{s:The basic framework and concepts}
%%%%%%%%%%%%%%%%%%%%%%%%%%%%%%%%%%%%%%%%%%%%%%%%%%%%%%%%%%%%%%%%

     Shelah in ~\cite{Sh 600} introduced  the axiomatic framework for the
notion of \emph{good frame};  his goal was to axiomatize
superstability. Below we offer a much simpler (and more general)
axiomatic setting we call \emph{weak forking} that in the
first-order case corresponds to simplicity.

\begin{definition}  A pair $\l\K,\nonfork \r$
is a \emph{weak forking notion}\index{weak forking} iff $\K$ is an AEC
and $\nonfork$ is a four-place relation called \emph{non-forking}
\index{non-forking}
$A\nonfork_{C}^N B$ for $C\subset A, B\subset N$,  $A,B,C\in \Ab(\K)$
such that $\nonfork$ satisfies
\begin{enumerate}
\item
\emph{Invariance}:
If $f:N\to N'$ is a $\K$-embedding, then
$A\nonfork_{C}^N B$ if and only if $f(A)\nonfork_{f(C)}^{f(N')}
f(B)$.

\item
\emph{Monotonicity:}
If $C\subset C' \subset B'\subset B$ and $N\prec N'$, then
$$
A\nonfork_{C}^N B \textrm{ if and only if } A\nonfork_{C'}^{N'} B'.
$$

\item \emph{Disjointness:}
\[
A\nonfork_{C}^N B \implies A\cap B\subseteq C.
\]

\item
\emph{Extension of independence:}
\index{Extension property}
If $M_0\prec M$ and $N_0\succeq M_0$, then there is a model $N\in \K$,
$M\prec N$,
and $f: N_0 \to N$ such that $f\restriction M_0 =\id _{M_0}$ and
$M\nonfork_{M_0}^{N} f(N_0) $.

\item
\emph{Continuity:}
If $\delta$ is a limit ordinal, $\{N_i\mid i<\delta\}$ is an increasing
continuous chain, and $M\nonfork _{M_0}^N N_{i}$ for $i<\delta$, then
$M\nonfork _{M_0} N_\delta$.

\item
\emph{Symmetry:}
\index{Symmetry}
if $M\nonfork_{M_0}^N N_0$, then $N_0 \nonfork_{M_0}^N M$.

\item
\emph{Transitivity:}
if $M\nonfork_{M_0}^N M_1$ and $M\nonfork_{M_1}^N M_2$, then $M
\nonfork_{M_0}^N M_2$.

\item
\emph{Local character:}
There is a cardinal $\kappa=\kappa(\K)$ such that
for any amalgamation base $A'\subset A\cup B$ there is an amalgamation
base $B'\subset B$, $|B'|=\kappa+|A'|$, with $A'\nonfork ^{N}_{B'} B$.

\item
\emph{Definability:}
There is a family $\mathcal{F}$ of $\kappa(\K)$ functions, definable
set-theoretically, such that
$A'\subset A\cup B$ is closed under $\mathcal{F}$,
then $A'\nonfork ^N _{B\cap A'} B$.

\end{enumerate}
\end{definition}

\begin{remark}
Axiom 9 is a very mild strengthening of the local character axiom.
It hides a brute force construction similar to the one in
Lemma~\ref{Skolem} and possible in the known examples. Suppose
local character holds, and that dependence relation makes sense for all sets.
Fix a well-ordering of the universe of $N\in \K$. For $\a\in A$, $\ell(\a)=n$,
define $\{f^n_i(\a) \mid i<\kappa(\K)\}$ to be an enumeration of the set
$B'\subset N$ such that $\a \nonfork ^{N}_{B'} N$. Letting
$\mathcal{F}:=\bigcup\{ f^n_i\mid i<\kappa,n<\omega\}$, we get the desired family.

The property stated in Axiom 9 was extracted from Section 4 of
Shelah's~\cite{Sh 87b}.

Of course, the local character property follows from definability of
independence.
\end{remark}

Axiom 9 and transitivity immediately give the following useful version
of the definability property.

\begin{claim}
There is a family $\mathcal{F}$ of $\kappa(\K)$ functions, definable
set-theoretically, such that
if $A\supset C$, $A\nonfork ^N _C B$, and $A'\subset A$ is closed
under $\mathcal{F}$, then $A'\nonfork ^N _{C\cap A'} B$.
\end{claim}

\begin{remark}
While we assume that the independence relation $\nonfork$ is defined over
amalgamation bases,
it is enough, for our purposes, to demand that the main properties of
independence such as
symmetry, transitivity, and extension holds only over models.

The extension property for the class follows from the amalgamation
assumptions we are making on
the class, see Section 2.
\end{remark}

\begin{remark}
To see that Shelah's notion of good frame is much more stronger than our,
imagine that $\K=\Mod(T)$ when $T$ is a complete first-order theory and
$\nonfork_{}$ is the usual first-order forking.  $\K$ is a good frame iff
$T$ is superstable, while $\l\K,\nonfork_{}\r$ is a weak forking notion iff
$T$ is simple.
\end{remark}

In the formulation of extension property, if $M_0=N_0$, we obtain existence
property of
independence. Let us state a form of the extension of independence property
that will be
useful later:

\begin{lemma}\label{useful ext}
If $M\nonfork_{M_0}^N N_0$ and $\bar N_0\succ N_0$, then there is a model
$\bar N\in \K$,
$N\prec \bar N$,
and $f:\bar N_0 \to \bar N$ such that $f\restriction N_0 =\id _{N_0}$ and
$M\nonfork_{M_0}^{\bar N} f(\bar N_0)$.
\end{lemma}

\begin{proof}
Applying extension of independence to $N$, $N_0$, and $\bar N_0$, we get a
model $\bar N\succ N$
and $f:\bar N_0\to \bar N$, identity over $N_0$, such that $N\nonfork
_{N_0}^{\bar N} f(\bar N_0)$.
Using symmetry and monotonicity we get $M\nonfork _{N_0}^{\bar N} f(\bar
N_0)$, and now symmetry
and transitivity give $M\nonfork_{M_0}^{\bar N} f(\bar N_0)$.
\end{proof}

\begin{examples}
\noindent
\begin{enumerate}
\item
Let $\K:=\Mod(T)$ when $T$ is a first-order complete theory, $\prec_{\K}$
is the usual
elementary submodel relation and
     $\nonfork_{}$ is the non-forking relation.   Clearly $\l\K,\prec_{\K}\r$
is
a weak forking notion iff $T$ is simple.  $\kappa$ in this case is
$\kappa(T)$.

It is not difficult to see that  $\l\K,\prec_{\K}\r$ is
a weak forking notion with $\kappa=\aleph_0$ iff $T$ is super-simple.

\item
Let $T$ be a countable first-order theory, and let
\[
\K^a:=\{M\models T \mid \ftp(\a/\emptyset,M) \text{ is an isolated type for
every }\a\in
|M|\}.
\]
A type $p\in S(A)$ is called \emph{atomic} iff   $A\cup\{\a\}$ is atomic
subset of $\C$ and $\a\models p$.

Suppose that $T$ is \emph{$\aleph_0$-atomically
stable}, i.e. for $\R[p]<\infinity$ for every atomic type, where
\begin{definition}
For $M\in \K^a$ and $\a\in M$ define by induction of $\alpha$ when
$\R[\phi(\x;\a)]\geq
\alpha$

$\alpha=0;\; M\models\exists\x\phi(\x;\a)$

For $\alpha=\beta+1$;

There are  $\b\supseteq \a$ and $\psi(\x;\b)$ such that
     \[
\R[\phi(\x;\a)\wedge\psi(\x;\b)]\geq\beta
\]
\[
\R[\phi(\x;\a)\wedge\neg\psi(\x;\b)]\geq\beta  \quad \text{and for every }
\c\supseteq \a
\]
there is $\chi(\x;\c) \text{ complete s.t. }$
\[
\R[\phi(\x;\a)\wedge\chi(\x;\c)]\geq\beta
\]
\end{definition}

An atomic set $A\subseteq \C$ is \emph{good}
iff for every consistent $\phi(\x;\a)$ (with $\a\in A$) there is
an isolated type $p\in S(A)$ containing $\phi(x;\a)$.
In the atomic case the countable good sets are amalgamation bases (compare
with Definition~\ref{def: AB}).  This follows from:

\begin{fact}[\cite{Sh 87a}]
Suppose  $A$ is countable. Then $A$ is good if and only if there is a
universal model over $A$.
\end{fact}

Suppose $A\cup B\cup C$ are inside $N\in \K^a$ and $C$ is good.
We let $A\nonfork _C B$ if for each $\a\in A$, $\ftp (\a/B)$
does not split over some finite subset of $C$.
Then $\l \K^a,\nonfork \r$ is a weak forking notion.

%\begin{notation}
%
%\[
%D_A:=\{\tp(\a/A)\mid A\cup \{\a\} \text{ is atomic}\}.
%\]
%\end{notation}
%
%\begin{fact}[\cite{Sh 87a}]
%If $|D_A|<2^{\aleph_0}$ for a countable $A$, then $A$ is good.
%\end{fact}

\item  Let
$\K$ be the class of elementary submodels of a totally transcendental
sequentially homogeneous model.
     Let $M_1\nonfork_{M_0}M_2$ stand for $\ftp(\a/M_2)$ does not
strongly-split over $M_0$ for every $\a\in |M_1|$.

Then $\l \K,\nonfork\r$ is a weak forking notion.
\end{enumerate}
\end{examples}

Compare the following with XII.2 of \cite{Sh c}.

\begin{definition}[Stable  systems]\label{st_sys}
Let $\l\K,\nonfork_{}\r$ be  weak forking notion.
     Suppose $I\subseteq \Pm(n)$, suppose ${\mathbf S}=\{M_s\mid s\in
I\cup\{n\}\}$
is a
$(\lambda,n)$-system.  The system $\mathbf  S$ is called
\emph{$(\lambda,I)$-stable in $M^{\S}_n$} if and only if
\begin{enumerate}

\item $A^{\mathbf S}_{s}$ is an amalgamation base for all $s\in I$,

\item
for all $s\in I$, for all $t\subseteq s$
\[
M^{\mathbf S}_{t}\;
\nonfork_{A^\S_{t}}^{M^{\mathbf S}_n}\;\;{\bigcup_{\substack {w\subseteq
s\\
w\not\supseteq t}}|M^{\mathbf S}_{w}|}.
\]
\end{enumerate}
\end{definition}

We make one more assumption on the $\l \K, \nonfork\r$.

\begin{axiom}[Generalized Symmetry]
     Let $\l\K,\nonfork_{}\r$ be  weak forking notion.  We say that
$\l\K,\nonfork_{}\r$
has the \emph{$(\lambda,n)$-symmetry property}  if a system ${\mathbf
S}=\{M_s\mid s\in
\P(n)\}$,
$\S\subset \K_\lambda$, is stable inside $M_n$ whenever there exists an
enumeration
$\bar s:=\l s(i)\mid i< 2^n-1 \r$ of $\Pm(n)$ (always without repetitions
such that $s(i_1)\subset s(i_2)\implies i_1<i_2$) such that
\begin{enumerate}
\item
$A^{\mathbf S}_{s(i)}$ is an amalgamation base for all $i$;

\item
\[
M^{\mathbf S}_{s(j)} \nonfork_{A^{\mathbf S}_{s(j)}}^{M^{\mathbf S}_n}
    {\bigcup_{i<j}|M^{\mathbf S}_{s(i)}|}.
\]
\end{enumerate}

\end{axiom}

In other words, under the generalized symmetry to get stability of the
$\Pm(n)$-system it is enough
to check the independence of just one ``face'' from the rest of the
$n$-dimensional cube, not all the faces as in the Definition~\ref{st_sys}.

We now state the generalized amalgamation properties, we omit the
superscripts $\S$ when the identity of the system is clear.

\begin{definition}[$n$-existence]
     Let $\l\K,\nonfork_{}\r$ be  weak forking notion. $\K$ has the
\emph{$(\lambda,n)$-existence property}
\index{$(\lambda,n)$-existence property}\index{existence property}
iff for every $(\lambda,\Pm(n))$-system
$\mathbf {S}=\l M_s\mid s\in \mathcal P^-(n)\r$
such that $\{M_t\mid t\subseteq s\}$ is a stable $(\lambda,|s|)$-system
for all $s\in \Pm(n)$,
there exists a model $M_n$ and $\K$-embeddings $\{f_s \mid s\in \Pm(n)\}$
such that
\begin{enumerate}
\item
$\{f_s(M_s) \mid s\in \Pm(n)\}\cup \{M_n\}$ is a stable system indexed by
$\P(n)$.
\item
the embeddings $f_s$ are coherent: $f_t\restriction M_s=f_s$ for
$s \subset t\in \Pm(n)$.
\end{enumerate}
\end{definition}

\begin{remark}
Let us clarify what is going on in the case $n=3$. We are given
the models $M_\emptyset$, $\{M_i\mid i<3\}$ and $\{M_{ij}\mid i<j<3\}$
Such that $M_i\nonfork ^{M_{ij}}_{M_\emptyset} M_j$ for all $i<j<3$.

The 3-existence property asserts that the three models can be embedded
into $M_{012}$ in a coherent way so that the images form a stable system
inside $M_{012}$. Note that this fails even in the first order case.

Failure of $(\aleph_0,3)$-existence is witnessed by the example
of a triangle-free random graph. Start with a triple of models $M_i$,
$i<3$ extending some $M_\emptyset$, and fix some elements $a_i\in M_i$.
Choose models $M_{01}$, $M_{02}$, and $M_{12}$ so that
$M_i\nonfork ^{M_{ij}}_{M_\emptyset} M_j$ for all $i<j<3$, and such
that
$M_{ij}\models R(a_i,a_j)$ for $i<j<3$. The system cannot be completed
since the model $M_{012}$ would witness a triangle.

This is an example of a non-simple first order
theory. It can be generalized to a failure of $(\aleph_0,n+1)$-amalgamation
by using $n$-dimensional tetrahedron-free graphs. Those examples are
simple first order theories.
\end{remark}

\begin{definition}[weak $n$-uniqueness]
Let $\S=\{M_s \mid s\in \Pm(n)\}$, $\S'=\{M'_s \mid s\in \Pm(n)\}$ be stable
systems of models in $\K$, where without loss of generality we assume
$M_\emptyset=M'_\emptyset$. We say that $\S$ and $\S'$ are \emph{piecewise
isomorphic}
if there are   $\{f_s:M_s \cong M'_s  \mid s\in \Pm(n)\}$,
where $f_\emptyset = \id _{M_{\emptyset}}$
and  $f_t\restriction M_s  = f_s$ for $s\subset t$.

Let $\l\K,\nonfork_{}\r$ be  weak forking notion. We say $\K$ has the
\emph{weak $(\lambda,n)$-uniqueness  property}
\index{weak $(\lambda,n)$-uniqueness}\index{weak uniqueness}
if the following holds. For any two $(\lambda,n)$-stable systems
$\mathbf {S}=\l M_s\mid s\in \mathcal P(n)\r$ and
$\mathbf {S}'=\l M'_s\mid s\in \mathcal P(n)\r$
such that $S\setminus \{M_n\}$ and $\S'\setminus \{M'_n\}$ are piecewise
isomorphic there are $M^*\in \K_\lambda$ and $\K$-embeddings $g:M_n\to M^*$
and $g':M'_n\to M^*$
such that $g(M_s)=g'(f_s(M_s))$ for all $s\in \Pm(n)$.
\end{definition}

\begin{remark}
In \cite{Sh 87b}, Shelah states a variant of weak $(\lambda,n)$-uniqueness
property. Shelah calls the property
\emph{failure of $(\lambda,n)$-non-uniqueness}, it is stated in item (2)
of Proposition~\ref{wu=failure of non uniq}. We show
that weak $(\lambda,n)$-uniqueness condition is equivalent to the
failure of $(\lambda,n)$-non-uniqueness.
\end{remark}

\begin{proposition}\label{wu=failure of non uniq}
Let $\l\K,\nonfork_{}\r$ be  weak forking notion. Then the following are
equivalent:
\begin{enumerate}
\item
$\K$ has the weak $(\lambda,n)$-uniqueness property;

\item
for every stable system
$\mathbf {S}=\l M_s\mid s\in \mathcal P^-(n)\r\subseteq \K_\lambda$ inside
some $M_n$ we have that $A^{\mathbf S}_n\in \Ab_\lambda(\K)$.
\end{enumerate}
\end{proposition}

\begin{proof}
If the weak uniqueness holds, then clearly the set $A_n$ is an amalgamation
base; we can take
the identity isomorphisms as the ``piecewise'' embeddings.

Now the converse. Let $\S^1$, $\S^2$ be piecewise isomorphic
stable systems indexed by $\Pm(n)$, inside $M^1_n$ and $M^2_n$
respectively. To show the weak uniqueness, it is enough to
construct a model $N^2_n$ and $g:M^1_n\isom N^2_n$ such that
$g\supset f_s$, $s\in \Pm(n)$ (it is enough to consider only the
$(n-1)$-element subsets $s$). Indeed, by invariance the system
$\S^2$ is stable inside $N^2_n$; by~(2) then
there are $M^*$ and $h_M: M^2_n\to M^*$, $h_N:N^2_n\to M^*$ over
$\S^2$. Then $h_N\circ g :M^1_n\to M^*$ and $h_M: M^2_n\to M^*$
are the needed embeddings.

The construction of $N^2_n$ and $g$ is a slight generalization of
the construction in the proof of Fact~\ref{pull-back}. As the
universe of $N^2_n$ we take the following set:
$$
|N^2_n|:= \Union_{s\in \Pm(n)} |M^2_s| \cup \left(|M^1_n|\setminus
\Union_{s\in \Pm(n)} |M^1_s|\right).
$$
Define the structure on the $|N^2_n|$ by copying it from the
structure $M^1_n$. Take a tuple $\a\in |N^2_n|$, it can be
uniquely presented as $\a=\cup_{s\in\Pm(n)} f_s(\a_s)\cup \b$,
where $\a_s\in M^1_s$ and $\b\in |M^1_n|\setminus \Union_{s\in
\Pm(n)} |M^1_s|$. For a relation $R\in L(\K)$ define
$$
N^2_n\models R(\a) \textrm{ if and only if } M^1_n\models
R(\cup_{s\in \Pm(n)}(\a_s)\cup\b).
$$
By construction $M^1_n$ is isomorphic to $N^2_n$.
\end{proof}

%\begin{definition}[$n$-uniqueness]
%   Let $\l\K,\nonfork_{}\r$ be  weak forking notion, it has the
%\emph{$(\lambda,n)$-uniqueness  property}
%\index{$(\lambda,n)$-uniqueness}\index{uniqueness}
%iff for every stable system
%$\mathbf {S}=\l M_s\mid s\in \Pm(n)\r$ the model $M_n$ is unique
%over $A^\S_n$.  COMPLETE ?????????????????????
%\end{definition}

\begin{definition}[goodness]
Let $\l\K,\nonfork_{}\r$ be  weak forking notion, it has the
\emph{$(\lambda,n)$-goodness
property}
\index{$(\lambda,n)$-goodness}\index{goodness}
iff
     $\l\K,\nonfork_{}\r$
has the {$(\lambda,n)$-symmetry property} and
%for every stable system
%$\mathbf {S}=\l M_s\mid s\in \mathcal \Pm(n)\r$ of models of cardinality
%$\lambda$,
has  the $(\lambda,n)$-existence property and the
weak $(\lambda,n)$-uniqueness  property.
\end{definition}

\begin{theorem}[characterizing goodness for f.o.]
Let $T$ be a complete countable first order  theory.  Suppose $T$ is
superstable
without dop
If
$\mathbf {S}=\l M_s\mid s\in \mathcal P^-(n)\r$ is a stable system
of models of cardinality
$\aleph_0$ then the following are equivalent:
\begin{enumerate}

\item
     the set $A^{\mathbf S}_n$ is an amalgamation base
\item
There is a prime and minimal model over $A^{\mathbf S}_n$.
\end{enumerate}
\end{theorem}

\begin{definition}[excellence]  Let $\l\K,\nonfork_{}\r$ be  weak forking
notion
and let $\lambda\geq\LS(\K)$.
$\l\K,\nonfork_{}\r$ is \emph{$\lambda$-excellent}
\index{$\lambda$-excellent}\index{excellent}
iff $\l\K,\nonfork_{}\r$  has the $(\lambda,n)$-goodness
property for every $n<\omega$.  When $\lambda=\LS(\K)$ we say that $\K$
excellent
     instead of
$\lambda$-excellent.
\end{definition}

\begin{theorem}[Shelah 1982]
Let $T$ be a complete countable first order  theory.  Suppose $T$ is
superstable
without DOP. Then the following are equivalent:
\begin{enumerate}

\item
     $\l\Mod(T),\prec\r$ is excellent.
\item
$\Mod(T)$ has the $(\aleph_0,2)$-goodness property.

\item  $T$ does not have the OTOP.
\end{enumerate}
\end{theorem}

For proof see  \cite{Sh c}.

\begin{fact}[Hart and Shelah 1986]
For every $n<\omega$
there is  an $\aleph_0$-atomically stable class $\K_n$ of atomic models of
a countable
f.o. theory such that $\K$ is has the $(\aleph_0,k)$-goodness property
for all $k<n$ but is not excellent.
\end{fact}

In section 3 we will prove that the existential quantifier in the
definition of excellent class can be
replaced with a universal quantifier:

\begin{theorem} If
$\l\K,\nonfork_{}\r$ is excellent then it has the
$(\lambda,n)$-goodness
property for every $n<\omega$ and every $\lambda\geq \LS(\K)$.
\end{theorem}

\begin{proof}
Immediate from Theorems \ref{ex transfer} and \ref{uq transfer}.
\end{proof}

%%% end section 1
%%% begin changed part:

%%%%%%%%%%%%%%%%%%%%%%%%%%%%%%%%%%%%%%%%%%%%%%%%%%%%%%%%%%%%%%%%%%%%
\section{A sufficient condition for Tameness}
%%%%%%%%%%%%%%%%%%%%%%%%%%%%%%%%%%%%%%%%%%%%%%%%%%%%%%%%%%%%%%%%%%%%
We start by explaining the main idea for obtaining
$(\lambda,\lambda^+)$-tameness from weak $(\lambda,2)$-uniqueness
and $(\lambda,2)$-existence.
We outline the general construction
and the induction step by a picture and later give a completely
formal argument.

Suppose $(a_1,M,N^1)\in p$ and $(a_2,M,N^2)\in q$ and
their restriction on small submodels of $M$ are equal.  Pick
$\{N^\ell_\alpha\prec_{\K}N_\ell\mid \alpha<\lambda\}\subseteq
\K_{<\lambda}$ increasing and continuous resolutions of $N_\ell$ and
$\{M^\alpha\prec_{\K}M\mid \alpha<\lambda\}\subseteq \K_{<\lambda}$
increasing resolution of $M$ such that $M_\alpha\prec
N^\ell_\alpha$, require that $a_\ell\in N^\ell_0$.

By the assumption there exist $N^*_0\succ_{\K}N^2_0$ of
cardinality less than $\lambda$ amalgam of $N^1_0$ and $N^2_0$ over $M_0$
mapping $a_1$ to $a_2$.

Our goal is to find models $\bar N^\ell\succ_{\K}N^\ell$ and
$\bar f_\lambda:\bar N^1\cong \bar N^2$ such that $\bar f_\lambda
(a_1)=a_2$.

The construction of the models and the mapping will be by induction on
$i<\lambda$ such that the following diagram commutes.

For $i=0$; let $\bar N_0^1\succ_{\K}N_0^1$ be an amalgam of $N_0^2$ and
$N_0^1$
over $M_0$ such that $g_0:N_0^2\rightarrow \bar N^1_0$, $g_0\restriction
M_0=\id_{M_0}$ and $g_0(a_2)=a_1$.  By  Fact~\ref{pull-back} there are
$\bar N_0^2\succ_{\K}N_0^2$ and  $\bar f_0:\bar
N_0^1\cong
\bar N_0^2$ such that $\bar f_0\supseteq g_0^{-1}$.  Using a strong form of
the
extension property (see Lemma~\ref{lambda+,lambda} below)
after renaming $\bar N_0^\ell$ we may assume that there exists $\hat
N^\ell_0\succ_{\K}\bar
N^\ell_0$ of cardinality $\lambda$ such that
$M\nonfork_{M_0}^{\hat N_0^\ell} \bar N^\ell_0$.
Using  Lemma \ref{lambda+,lambda} once more we find $\check
N_1^\ell\in \K_\lambda$ such that
$\hat N^\ell_0 \nonfork_{M_0}^{\check N_1^\ell}N^\ell_1$. Since
$\check N^\ell_0\succ \hat N^\ell_0$, by monotonicity
$$
(*)\qquad M\nonfork_{M_0}^{\check N_0^\ell} \bar N^\ell_0.
$$
Now take $\tilde N^\ell_1\prec \check N^\ell_1$ of
cardinality $\lambda$ such that it contains
$ |\bar N^\ell_0|\cup |N^\ell_{1}|$.

Monotonicity applied to $(*)$ gives that  $M_1\nonfork_{M_0}^{\tilde
N^\ell_1}\bar N^\ell_0$ holds for $\ell=1,2$. An application of the weak
2-uniqueness property produces a model $\bar N^1_1\succ \tilde N^1_1$ and
$g_1:\tilde N^2_1\rightarrow \bar N^1_1$ such that $g_1\supseteq \bar
f_0^{-1}$
and $g_1\restriction M_1=\id_{|M_1|}$.  Now using
Fact
\ref{pull-back} there are $\bar N_1^2\succ_{\K}\tilde N_1^2$ and  $\bar
f_1:\bar
N_1^2\cong
\bar N_1^1$ such that $\bar f_1\supseteq g_1^{-1}$.

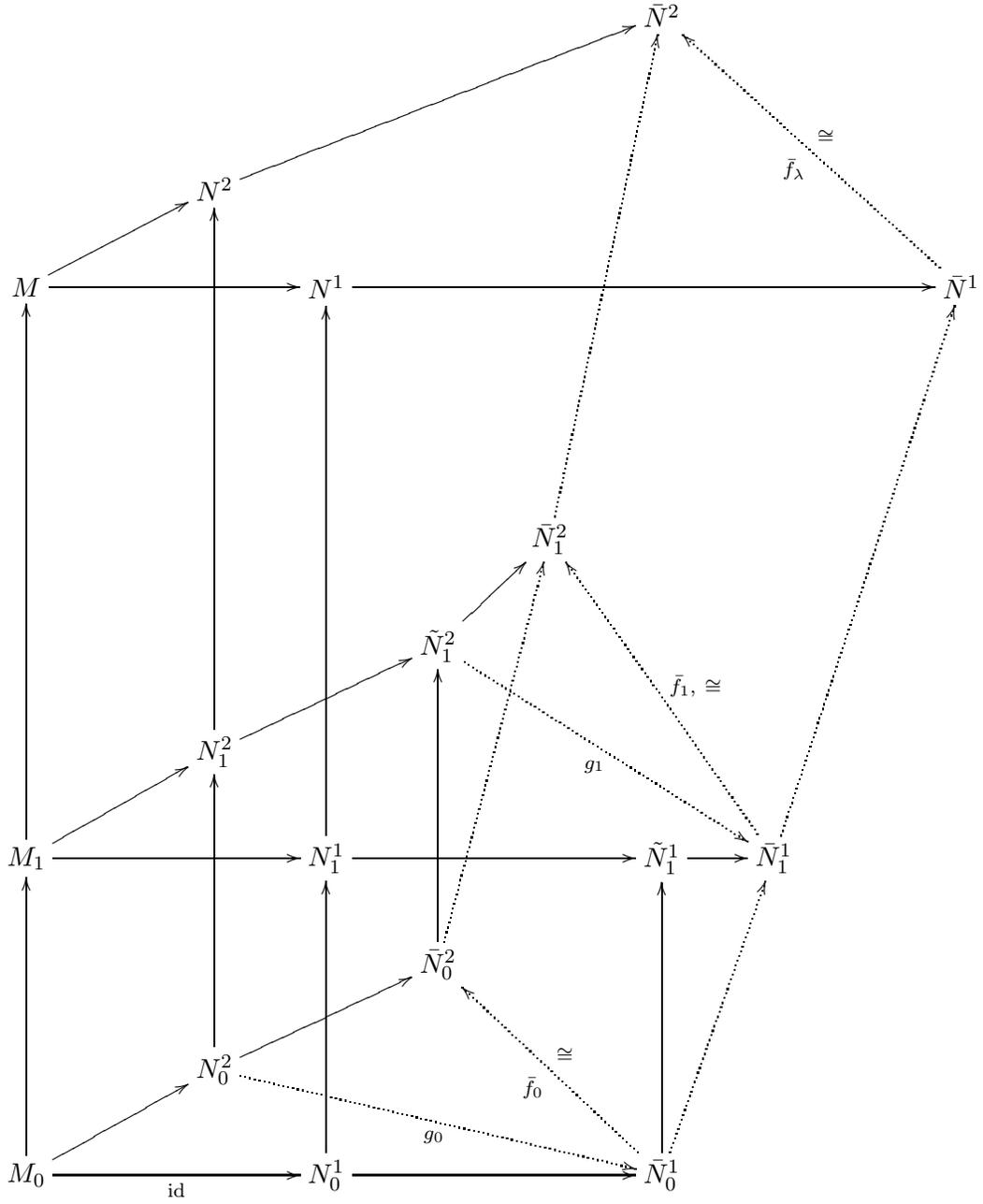
\[
     \xymatrix{\ar @{}  &&&&&& {\bar N}^2    &&&&
\\ \\
&& N^2  \ar[rrrruu]^{}
\\
M \ar[rrr] _{} \ar[rru]^{}  & && N^1 \ar[rrrrrr] _{} & & & &&&{\bar N}^1
     \ar@{.>}[uuulll]^{\bar f_\lambda}_{\cong}
\\  \\ \\
&&&&& {\bar N}^2_1 \ar@{.>}[uuuuuur]&&&&
\\
&&&& {\tilde N}^2_1 \ar@{.>}[rrrdd]_{g_1} \ar[ur]^{} &&&&&
\\
&& N_1^2 \ar[rru] \ar[uuuuuu]_{}
\\
     M_1 \ar[uuuuuu]^{} \ar[rrr] _{} \ar[rru]^{}  & && N_1^1
\ar[uuuuuu]_{} \ar[rrr] _{}
     & & &{\tilde N} ^1_1
%  \ar[uuul]^{f_1}
\ar[r]_{} & {\bar N} ^1_1 \ar@{.>}[uuuuuurr]^{}\ar@{.>}[uuull]_{\bar f_1,\
\cong}
\\
&&&& {\bar N}^{2}_0 \ar@{.>}[uuuur]^{}\ar[uuu] &&&&&
\\
&& N_0^2 \ar@{.>}[rrrrd]_{g_0}\ar[rru]^{}   \ar[uuu]_{}
\\
M_0\ar[uuu]^{} \ar[rru]^{}  \ar[rrr]_{\id} & && N_0^1
\ar[uuu]_{} \ar[rrr]_{} & & & {\bar N}_0^{1} \ar[uuu]_{}
\ar@{.>}[uull]^{\bar f_0}_{\cong} \ar@{.>}[uuur]_{}
     }
\]

%\[
%  \xymatrix{\ar @{}  [dr] &&&& {\bar N}^2    &&&&& \\
%  && N^2  \ar[rru]^{\id} \\
%  M
%  \ar[rrr] _{\id} \ar[rru]^{\id}  & && N^1 \ar[rrrr] _{\id} & & & & {\bar
%N}^1
%  \ar[uulll]^{\bar f_\lambda}_{\cong}
%   \\  \\ \\
%\\ &&&& {\bar N}^2_1 \ar[uuuuuu]^{\id}  &&&&& \\
%&& N_1^2 \ar[rru] \ar[uuuuuu]_{\id}\\
%  M_1 \ar[uuuuuu]^{\id} \ar[rrr] _{\id} \ar[rru]^{\id}  & && N_1^1
%\ar[uuuuuu]_{\id} \ar[rrrr]
%_{\id}
%  & & && {\bar N} ^1_1\ar[uuuuuu]^{\id}
%  \ar[uulll]^{\bar f_1}_{\cong}
%\\ &&&& {\bar N}^{2}_0 \ar[uuu]^{\id} &&&&& \\
%&& N_0^2 \ar[rru]^{}   \ar[uuu]_{\id} \\
%M_0\ar[uuu]^{\id} \ar[rru]^{}  \ar[rrr]_{\id} & && N_0^1 \ar[uuu]_{\id}
%\ar[rrrr]_{\id}
%& & & & {\bar N}_0^{1} \ar[uuu]_{\id}   \ar[uulll]^{\bar f_0}_{\cong}
%  }
%\]

For the rest of this section we deal with $(\lambda,2)$-existence and weak
uniqueness
properties.
We make it explicit that
$(\lambda,2)$-existence is
simply the $\lambda$-extension property for independence (see
Definition~\ref{lambda-ext});
and weak $(\lambda,2)$-uniqueness corresponds to the first-order
stationarity.
This makes transparent the argument showing, for example,
$(\lambda^+,\lambda)$-tameness from
$(\lambda,2)$-existence and weak uniqueness. The first-order relativization
of the proof
goes along these lines: let $p$, $q$ be types over $M$ of size $\lambda^+$
that agree over
all $\lambda$-submodels of $M$. With $\lambda\ge \kappa(\K)$, by local
character we can find
$M_0\prec M$, $\|M_0\|=\lambda$, such that $p$, $q$ do not fork over $M_0$.
By assumption
$p\restriction M_0=q\restriction M_0$, so stationarity gives $p=q$. Of
course this outline avoids
several important issues; for example, we assume ``stationarity'' only in
$\lambda$, and we used
$\lambda^+$-stationarity in the argument above.

Let us restate the definitions of 2-existence and weak uniqueness here.

\begin{definition}
We say that $\K$ has \emph{$(\lambda,2)$-existence property} if for any
triple of models
$M_0\prec M_1,M_2$, all in $\K_\lambda$, there is a model $M^*\in
\K_\lambda$ and
$\K$-embeddings $f_\ell: M_i\to M^*$, $i\ell=1,2$,
such that $f_1\restriction M_0=f_2\restriction M_0$ and
$f_1(M_1) \nonfork^{M^*}_{f_1(M_0)} f_2(M_2)$.
\end{definition}

\begin{remark}
Equivalently, $(\lambda,2)$-existence property holds if there is $M^*\succ
M_1$ and a map
$f:M_2\to M^*$ over $M_0$ such that $M_1 \nonfork ^{M^*}_{M_0} f(M_2)$. So
$(\lambda,2)$-existence is really the $\lambda$-extension property:

\begin{definition}\label{lambda-ext}
If in the extension for independence property all the models are in
$\K_\lambda$, we say that $\K$ has the \emph{$\lambda$-extension property
for independence}.
\end{definition}
\end{remark}

\begin{definition}\label{def: 2-wu}
The class $\K$ has \emph{weak $(\lambda,2)$-uniqueness property} if for any
two $\lambda$-systems
$\S^\ell=\{M^\ell_0,M^\ell_1,M^\ell_2,M^\ell_3\}$, $\ell=0,1$, that are
stable (i.e.,
$M^\ell_1\nonfork^{M^\ell_3}_{M^\ell_0} M^\ell_2$) and piecewise isomorphic
(i.e., there are
$f_i:M^0_i\cong M^1_i$ for $i=0,1,2$ with $f_1\restriction
M_0=f_2\restriction M_0$), let
$f_0(M_0):=f_1(M_0)$ there is a
model $M^*\in \K_\lambda$ and embeddings $g^\ell:M^\ell_3\to M^*$ such that
$g^0(M_i)=g^1(f_i(M_i))$ for $i=0,1,2$.
\end{definition}

\begin{remark}\label{uq char}
(1) Equivalently, in $(\lambda,2)$-uniqueness we may demand that $M^*\succ
M^1_2$, i.e., $f^1$ is the identity
embedding.

(2) In the first-order case, weak 2-uniqueness says that ${\mathrm
tp}(M_1\cup
M_2/M_0)$ is uniquely determined
by ${\mathrm tp}(M_1/M_0)\cup {\mathrm tp}(M_2/M_0)$ as long as
$M_1\nonfork _{M_0} M_2$. So it really is the
analog of stationarity.

(3) Using the isomorphism axioms in the definition of AEC, weak
$2$-uniqueness can be viewed as an amalgamation property for 2
isomorphisms.

(4) The property is called \emph{weak} uniqueness since the property
$M_1\nonfork^N_{M_0}M_2$ (all models of cardinality $\lambda$) implies
existence of $N'\prec N$ such that $|N'|\supseteq |M_1|\cup |M_2|$ and
$N'$ is prime and minimal over the set $|M_1|\cup |M_2|$.  This stronger
property occurs in two different situations: Superstable countable
elementary classes without the OTOP (see Chapter XII in \cite{Sh c}) and
also in excellent classes of atomic models of a first-order theory (in
the sense of
\cite{Sh 87b}).

(5) Compare Definition \ref{def: 2-wu} with Definition 6.2 and
Claim 6.7 from \cite{Sh 576}.

\end{remark}

\begin{lemma}
Suppose $(\lambda,2)$-uniqueness holds. Let $\S^\ell$, $\ell=1,2$, be stable
and piecewise isomorphic $\lambda$-systems,
$\S^\ell=\{M^\ell_0,M^\ell_1,M^\ell_2,M^\ell_3\}$.
Then there are $\K_\lambda$-models $N^\ell\succ M^\ell_3$, $\ell=0,1$, and
$\bar f:N^1\cong N^2$ that extends the isomorphisms $f_i:M^1_i\to M^2_i$,
$i=0,1,2$.
\end{lemma}

\begin{proof}
Let $N^2\succ M^2_3$ and $f:M^1_3 \to N^2$ be as guaranteed by the weak
$(\lambda,2)$-uniqueness, in the sense of Remark~\ref{uq char}(1). Using
Fact~\ref{pull-back},
we get the needed $N^1\succ M^1_3$ and the isomorphism
$\bar f$ extending $f$, and therefore
all the mappings $f_i$, $i=0,1,2$.
\end{proof}

\begin{lemma}\label{lambda+,lambda}
Let $\lambda\ge \LS(\K)+\kappa(\K)$, and suppose $\lambda$-extension
property holds.
Let $\|M\|=\lambda ^+$, $M_0\prec M,N_0$, where
$M_0,N_0\in \K_\lambda$. Then there is $N\succ M$ and $f:N_0\to N$ over
$M_0$ such that
$M\nonfork _{M_0}^N f(N_0)$.
\end{lemma}

In short, $\lambda$-extension implies extension when one of the models has
size $\lambda^+$.

\begin{proof}
Let $\{M_i\mid i<\lambda^+\}$ be an increasing continuous chain of models
with
$\Union_{i<\lambda^+} M_i=M$, $\|M_i\|=\lambda$, and $M_0$ given in the
statement of the lemma.

By induction on $i<\lambda^+$, we build models $N_i$, $\|N_i\|=\lambda$ and
$\K$-embeddings
$f_{ij}: N_i \to N_j$ such that:
\begin{enumerate}
\item
$N_i\succ M_i$ for all $i<\lambda^+$;
\item
$\{N_i, f_{ij}\}$ form a directed system;
\item
$f_{ij}\restriction M_i = \id_{M_i}$;
\item
$M_i\nonfork^{N_i}_{M_0} f_{0i}(N_0)$.
\end{enumerate}
This is clearly sufficient: letting $N$ be the direct limit of
$\{N_i,f_{ij} \mid i<\lambda^+\}$,
we have $M\prec N$ by (1) and (3) and letting $f:=f_{0\lambda^+}$ we have
$M\nonfork _{M_0}^N f(N_0)$ and $f\restriction M_0=\id_{M_0}$ by (3) and
(4).

Now the construction: $N_0$ is given; having constructed $N_i$ and $f_{jk}$
for $j\le k\le i$
satisfying (1)--(4), build $N_{i+1}$ and $f_{j,i+1}$.

By $\lambda$-extension applied to $M_{i+1}\succ M_i$  and $N_i\succ M_i$,
there is a model $N_{i+1}\succ M_{i+1}$ and embedding $f_{i,i+1}:N_i\to
N_{i+1}$ such that
$f_{i,i+1}\restriction M_i =\id_{M_i}$ and $M_{i+1}\nonfork
^{N_{i+1}}_{M_i} f_{i,i+1}(N_i)$.
For $j<i$ we define $f_{j,i+1}:=f_{i,i+1}\circ f_{ji}$, and $f_{i+1,i+1}$
is the identity.
Thus, we have met (1)--(3).

We prove that we have (4). Since $M_{i+1}\nonfork ^{N_{i+1}}_{M_i}
f_{i,i+1}(N_i)$ and
$f_{0,i+1}(N_0)\prec f_{i,i+1}(N_i)$ by monotonicity we have
$M_{i+1}\nonfork ^{N_{i+1}}_{M_i} f_{0,i+1}(N_0)$. By induction hypothesis
$M_i\nonfork^{N_i}_{M_0} f_{0i}(N_0)$, so by invariance (applying
$f_{i,i+1}$) and monotonicity
$M_i\nonfork^{N_{i+1}}_{M_0} f_{0,i+1}(N_0)$. Symmetry and transitivity now
give the desired
$M_{i+1}\nonfork^{N_{i+1}}_{M_0} f_{0,i+1}(N_0)$.

Suppose now that $i$ is a limit ordinal. Let $N_i$ be the direct limit of
the system
$\{N_j, f_{jk}\mid j\le k<i\}$. By induction hypothesis
$M_j\nonfork^{N_j}_{M_0} f_{0j}(N_0)$ for all $j<i$. Applying $f_{ji}$ and
noting
$f_{ji}(f_{0j}(N_0))=f_{0i}(N_0)$, by invariance and monotonicity we get
$M_j \nonfork^{N_i}_{M_0} f_{0,i}(N_0)$ for all $j<i$. By continuity of
independence we finally have
$M_i\nonfork^{N_i}_{M_0} f_{0,i}(N_0)$.
\end{proof}

\begin{remark}\label{wma agreement}
The proof is actually a five-line argument if we phrase its key element
this way:

Given $N_i$ such that $M_i\nonfork _{M_0}^{N_i} N_0$,
by $\lambda$-extension property, \emph{we may assume} that
there is $N_{i+1}\succ N_i, M_{i+1}$ such that
$M_{i+1} \nonfork _{M_i}^{N_{i+1}} N_i$. By monotonicity $M_{i+1} \nonfork
_{M_i}^{N_{i+1}} N_0$,
and since also $M_i\nonfork _{M_0}^{N_{i+1}} N_0$, symmetry and
transitivity give
$M_{i+1}\nonfork _{M_0}^{N_{i+1}} N_0$.

So below we agree to use an appropriate ``we may assume'' in the place of a
directed
system argument. This makes the proofs much more transparent and does not
limit the
generality.
\end{remark}

\begin{corollary}
Let $\chi \ge \LS(\K)+\kappa(\K)$, and suppose $\mu$-extension property
holds
for all $\chi \le \mu< \lambda$. Let $\|M\|=\lambda $, $M_0\prec M,N_0$,
where
$M_0,N_0\in \K_\chi$. Then there is $N\succ M$ and $f:N_0\to N$ over $M_0$
such that
$M\nonfork _{M_0}^N f(N_0)$.
\end{corollary}

\begin{proof}
The same argument as in Lemma~\ref{lambda+,lambda}; the only difference is
that the
sequence $\{N_i\mid i<\lambda\}$ is such that $\|N_i\|=\chi+|i|$.
\end{proof}

\begin{theorem}\label{2.w.u.+ extension implies tame}
Suppose that $\K$ is an AEC with a weak forking notion.
Suppose for some $\chi \ge \LS(\K)+\kappa(\K)$ for all $\mu \in
[\chi,\lambda)$
weak $(\mu,2)$-uniqueness and $\mu$-extension hold. Then $\K$ is
$(\chi,\lambda)$-tame.
\end{theorem}

\begin{corollary}
Suppose that $\K$ is an AEC with a weak forking notion.
Suppose for some $\lambda \ge \LS(\K)+\kappa(\K)$
weak $(\lambda,2)$-uniqueness and $\lambda$-extension hold. Then $\K$ is
$(\lambda,\lambda^+)$-tame.
\end{corollary}

\begin{proof}[Proof of the theorem.]
Let $M\in \K$ be of size $\lambda$, and let $a_2$, $a_1$ have the same
Galois type over
every $\K_\chi$-submodel of $M$. We are constructing models $\bar N^\ell$
extending $N^\ell$
and a $\K$-isomorphism $\bar f:\bar N^2\to \bar N^1$ such that
$\bar f(a_2)=a_1$ and
$\bar f \restriction M=\id_M$.

Since this is the first time we are using our agreement from
Remark~\ref{wma agreement}, let us note that, strictly speaking,
the models $\bar N^\ell$, $\ell=1,2$,
arise as certain direct limits, $N^\ell$ embed into
$\bar N^\ell$ via $f^\ell$, and the condition is
$\bar f(f^2(a_2))=f^1(a_1)$.

Let $\{M_i \mid i<\lambda \}$ and $\{N^\ell_i \mid i<\lambda\}$,
$\ell=1,2$ be increasing continuous chains such that
\begin{enumerate}
\item
$M_i\prec M$ for all $i<\lambda$ and $\Union_{i<\lambda} M_i=M$;
\item
$\|M_i\|=\chi +|i|$ for all $i<\lambda$;
\item
$N^\ell_i\prec N^\ell$ for all $i<\lambda$ and $\Union_{i<\lambda}
N^\ell_i=N^\ell$;
\item
$M_i\prec N^\ell_i$ for all $i<\lambda$ and $\|N^\ell_i\|=\chi +|i|$;
\item
$a_\ell \in N^\ell_0$ and $M\nonfork ^{N^\ell}_{M_0} N^\ell_0$.
\end{enumerate}
This is easy since by local character we can find $N^\ell_0\prec N^\ell$
containing $a_\ell$ and $M_0\prec M$, $M_0\prec N^\ell_0$, such that
$\|M_0\|=\|N^\ell_0\|=\chi$ and $M\nonfork ^{N^\ell}_{M_0} N^\ell_0$.
The rest is immediate.

By induction on $i<\lambda$ we build increasing continuous chains
$\{\bar N^\ell_i\mid i<\lambda\}\subseteq \K_{<\lambda}$ and $\{\hat
N^\ell_i\mid
i<\lambda\}\subseteq \K_\lambda$ as well as isomorphisms $\bar f_i:\bar
N^2_i
\isom
\bar N^1_i$ such that
\begin{enumerate}
\item
$N^\ell_i\prec \bar N^\ell_{i} \prec \hat N^\ell_i$ for all $i<\lambda$,
$\ell=1,2$;
\item
$\|N^\ell_i\|=\|\bar N^\ell_{i}\|=\chi +|i|$ and $\|\hat N^\ell_i\|=\lambda$
for all $i<\lambda$;
\item
$\bar f_0(a_2)=a_1$ and $\bar f_i\subset \bar f_j$ for $i<j<\lambda$;
\item
$\bar f_{i}\restriction M_{i}=\id_{M_{i}}$ for all $i<\lambda$;
\item
$M \nonfork^{\hat N^\ell_i}_{M_i} \bar N^\ell_i$.
\end{enumerate}

Begin with $i=0$. Since $\gatp(a_2/M_0)=\gatp(a_1/M_0)$, there is a model
$\bar N^1_0 \succ N^1_0$
and an embedding $g_0:N^2_0 \to \bar N^1_0$.
Let $\bar N^2_0\in \K$ be such that $\bar N^2_0\succ N^2_0$
and $\bar N^2_0$ is isomorphic to $\bar N^1_0$ via some $\bar f_0$ such
that
$\bar f_0\restriction N^2_0 = g_0^{-1}$ (possible by Fact~\ref{pull-back}).

By extension, we may assume that there are $\hat N^\ell_0\succ M, \bar
N^\ell_0$
such that $M \nonfork _{M_0}^{\hat N^\ell_0} \bar N^\ell_0$.

For $\alpha$ a limit ordinal, let
$\bar N^\ell_\alpha:=\Union_{i<\alpha} \bar N^\ell_i$,
$\hat N^\ell_\alpha:=\Union_{i<\alpha} \hat N^\ell_i$,
and $\bar f_\alpha:=\Union_{i<\alpha} \bar f_i$.
It is routine to check that (1)--(4) hold, and we need to establish (5). By
the induction hypothesis
and monotonicity, for all $i<\alpha$ we have $M \nonfork^{\hat
N^\ell_\alpha}_{M_\alpha} \bar N^\ell_i$.
So by continuity we get $M \nonfork^{\hat N^\ell_\alpha}_{M_\alpha} \bar
N^\ell_\alpha$.

For the successor case, let $\mu:=\chi +|i|$.
%step 1
Since $N^\ell_i\prec N^\ell_{i+1}, \hat N^\ell_i$
by $(<\lambda)$-extension and Lemma~\ref{lambda+,lambda},
we can find $\check N^\ell_{i+1}$ of cardinality $\lambda$
such that $\hat N^\ell_i\nonfork^{\check N^\ell_{i+1}}_{N^\ell_i}
N^\ell_{i+1}$.
(Of course $N^\ell_{i+1}$ embeds into $\check N^\ell_{i+1}$, and
as in Remark~\ref{wma agreement} we assume the embedding is identity.)

Let $\tilde N^\ell_{i+1}\prec \check N^\ell_{i+1}$  be of cardinality $\mu$
such that
$|\tilde N^\ell_{i+1}|\supseteq
|\bar N^\ell_i|\cup |N^\ell_{i+1}|$.
By monotonicity, we still have
$M_{i+1}\nonfork^{\tilde N^\ell_{i+1}}_{M_i} \bar N^\ell_i$.
%step 2
By weak $(\mu,2)$-uniqueness (the systems $\{M_{i+1},M_i,\bar N^2_i\}$ and
$\{M_{i+1},M_i,\bar N^1_i\}$
are piecewise isomorphic inside $\tilde N^2_{i+1}$ and $\tilde N^1_{i+1}$),
there is a model $\bar N^1_{i+1} \succ \tilde N^1_{i+1}$
and an embedding $g_{i+1}: \tilde N^2_{i+1} \to \bar N^1_{i+1}$ that
extends the
identity map on $M_{i+1}$ and the isomorphism $\bar f_i$.
%step 3
Using Fact~\ref{pull-back} again, we get $\bar N^2_{i+1}\in \K_\mu$ such
that $\bar N^2_{i+1}\succ \tilde N^2_{i+1}$
and $\bar N^2_{i+1}$ is isomorphic to $\bar N^1_{i+1}$ via some
$\bar f_{i+1}$ such that
$\bar f_{i+1}\restriction \tilde N^2_{i+1} = g_{i+1}^{-1}$.
%step 4
By $(<\lambda)$-extension and
Lemma~\ref{lambda+,lambda}, we may assume that
there are $\hat N^\ell_{i+1}\succ \hat N^\ell_i, \bar N^1_{i+1}$ such that
$M \nonfork _{M_{i+1}}^{\hat N^\ell_{i+1}} \bar N^\ell_{i+1}$ for
$\ell=1,2$.

Having finished the construction, it remains to note that
$a_\ell\in \bar N^\ell:=\bar N^\ell_\lambda$, $N^\ell\prec \bar N^\ell$,
$\ell=1,2$, and the isomorphism $\bar f_\lambda:\bar N^2 \isom \bar N^1$
fixes $M$ and sends $a_2$ to $a_1$. Thus $\gatp(a_2/M)=\gatp(a_1/M)$.
\end{proof}

The following is a variation on Definition 0.23 from \cite{Sh 576}:

\begin{definition}
Let $\mu>\LS(\K)$.  The class $\K$ is called \emph{$\mu$-local} iff for
every $M\in\K_\mu$ and every resolution $\{M_i\prec_{\K}M \mid
i<\mu\}\subseteq\K_{<\mu}$ we have that
\[
(\forall i<\mu)[p\restriction M_i=q\restriction M_i\implies p=q]
\text{ for all }p,q\in\gaS(M).
\]
\end{definition}

It is easy to see that if an AEC is $\lambda^+$-local then it is
$(\lambda,\lambda^+)$-tame.  Notice that the proof of Theorem
\ref{2.w.u.+ extension implies tame} gives us the slightly
stronger result:

\begin{corollary}
Suppose that $\K$ is an AEC with a weak forking notion.
Suppose for some $\lambda \ge \LS(\K)+\kappa(\K)$
weak $(\lambda,2)$-uniqueness and $\lambda$-extension hold. Then $\K$ is
$\lambda^+$-local.
\end{corollary}

%%%end section 2

%%%%%%%%%%%%%%%%%%%%%%%%%%%%%%%%%%%%%%%%%%%%%%%%%%%%%%%%%%%%%%%%%%%%
\section{Stepping up}
%%%%%%%%%%%%%%%%%%%%%%%%%%%%%%%%%%%%%%%%%%%%%%%%%%%%%%%%%%%%%%%%%%%%

For this section, $\K$ is $\chi$-excellent; $\chi\ge \LS(\K)+\kappa(\K)$.
Our goal is to show that a $\chi$-excellent AEC $\K$ is
$(\chi,\infty)$-tame. For this,
it is enough to establish that

\begin{quote}
     In $\K$ $(\lambda,2)$-existence and weak uniqueness hold for $\lambda\ge
\chi$.
\end{quote}

This will follow from two theorems:

\begin{theorem}\label{ex transfer}
Suppose that $\lambda>\chi$ and $\K$ has
$(\mu,\le n+1)$-existence and weak $(\mu,n)$-uniqueness for all $\chi\le
\mu<\lambda$.
Then $\K$ has $(\lambda,\le n)$-existence.
\end{theorem}

\begin{theorem}\label{uq transfer}
Suppose that $\lambda>\chi$ and $\K$ has weak
$(<\lambda,\le n+1)$-uniqueness. Then $\K$ has weak $(\lambda,\le
n)$-uniqueness.
\end{theorem}

Let us start with some definitions and preliminary results.

\begin{definition}
Let $\S_\ell=\{M^s_\ell \mid s\in \Pm(n)\}$, $\ell=1,2$ be
$\Pm(n)$-systems such that $M^s_1\prec M^s_2$ for all $s\in \Pm(n)$.
We then write $\S_1\prec \S_2$.

If in addition for all $s\in\Pm(n)$, $|s|=n-1$, the $\S_1\cup
\S_2$-submodels
of $M^s_2$ form a stable $(\lambda,n)$-system inside $M^s_2$, then we say
that
$\S_1\prec \S_2$ are \emph{independent} and write $\S_1\nonfork \S_2$.
\end{definition}

Let us illustrate what the definition of $\S_1 \nonfork \S_2$
says in the simplest case when the dimension is 2, so
$\S_\ell=\{M^0_\ell,M^1_\ell,M^2_\ell,M^3_\ell\}$, $\ell=1,2$.
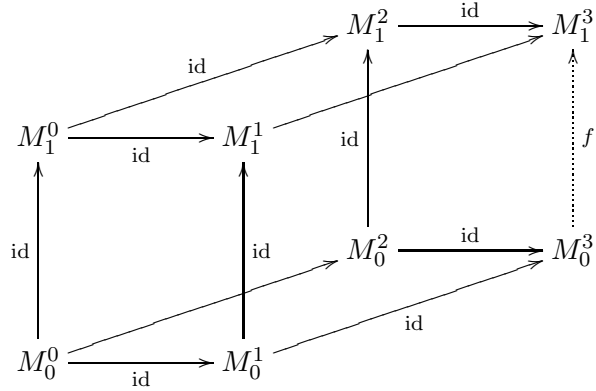
\[
     \xymatrix{\ar @{}
   &&& M^2_1 \ar[rr]^{\id} &&
M^3_1
\\
M^0_1\ar[rrru]^{\id}  \ar[rr]_{\id} &
   & M^1_1 \ar[rrru]
\\ &&& M^{2}_0 \ar[uu]^{\id} \ar[rr]^{\id} && M^3_0
\ar@{.>}[uu]_{f}
\\
M^0_0\ar[uu]^{\id} \ar[rrru]^{}  \ar[rr]_{\id} &
&  M_0^{1} \ar[uu]_{\id}  \ar[rrru]_{\id}
     }
\]
If $\S_1\nonfork \S_2$ then in $M^1_1$ we have
$M^0_1\nonfork^{M^1_1}_{M_0^0} M^1_0$, and similar for $M_1^2$.
However, this is not a 3-dimensional stable system in $M^3_1$ yet.
Existence of embedding $f$ that makes the system stable is
obtained in Lemma~\ref{cl1}
below. This is really a generalized extension property.

\begin{lemma}\label{cl1}
Let $\lambda\ge \LS(\K)+\kappa(\K)$, and suppose $(\lambda,\le
n+1)$-existence
and weak $(\lambda,n)$-uniqueness hold, $n\ge 2$.
Let $\S_1\prec \S_2$ be independent stable $(\lambda,n)$-systems inside the
models $M^{n}_1$, $M^{n}_2$
respectively. Then there is $\bar M^{n}_2\succ M^{n}_2$ and an embedding
$f:M^n_1 \to \bar M^n_2$
such that $f\restriction \Union _{M\in \S_1} M = \id$ and the system
$\S_1\cup \{f(M^n_1)\}\cup \S_2$
is a stable $(\lambda,n+1)$-system inside $\bar M^n_2$.
\end{lemma}

\begin{proof}
By $(\lambda,n+1)$-existence, there is $\bar M^n_2$ and embeddings
$f_s:M^s_2 \to \bar M^n_2$,
$s\in \Pm(n)$, $|s|=n-1$, and $f:M^n_1\to \bar M^n_2$ such that the images
form a stable
$(\lambda,n+1)$-system in $\bar M^n_2$.

Now in particular the image of $\S_2$ is a stable $(\lambda,n)$-system in
$\bar M^n_2$. So
by weak $(\lambda,n)$-uniqueness, the models $M^n_2$ and $\bar M^n_2$ can
be amalgamated
over $\S_2$. Thus by Fact~\ref{pull-back} we may assume that actually $\bar
M^n_2\succ M^n_2$.
Finally, $\bar M^n_2$ and $f$ are as needed.
\end{proof}

\begin{lemma}\label{resolution}
Let $\chi\ge \LS(\K)+\kappa(\K)$, $\lambda >\chi$. Let $\S=\{M^s\mid s\in
\Pm(n)\}$
be a stable $(\lambda,\Pm(n))$-system inside some $M^n$. There is a sequence
$\S_i=\{M^s_i \mid s\in \Pm(n)\}$, for $i<\lambda$ such that
\begin{enumerate}
\item
$\S_i$ is a $(\chi+|i|,\Pm(n))$-system;
\item
$\S_i \prec \S_{i+1}$ and $\S_i\nonfork \S_{i+1}$ for $i<\lambda$;
\end{enumerate}
\end{lemma}

\begin{proof}
Let $\chi_0$ be large enough regular so that $H(\chi_0)$ contains all the
information about the system
$\S$. Let $\B_i\prec \l H(\chi_0),\in \dots \r$ be an internal chain of
models, with
$\|\B_i\|=\chi+|i|$, and such that $(M^\emptyset)^{\B_i}$ has size
$\chi+|i|$.

By definability of independence, $\S_i:=\S^{\B_i}$ is a stable
$\Pm(n)$-system. It remains to show (2).
Let $s\in \Pm(n)$, $|s|=n-1$, let $j:=i+1$, and let $\mu:=\chi+|i|$. We are
showing that
$\{M^t_i \mid t\subseteq s\}\cup \{M^t_j \mid t\subsetneq s\}$ is a stable
$(\mu,n+1)$-system in $M^s_j$.

By generalized symmetry, it is enough to show that $M^s_i \nonfork
^{M^s_j}_{A^s_i} A^s_j$.
But this follows from definability of independence: $A^s_i=A^s_j\cap M^s_i$
and $M^s_i$ is closed
under the $\kappa$-many functions that define independence.
\end{proof}

\begin{remark}
This is the only place where we had to use the generalized symmetry axiom.
\end{remark}

\begin{lemma}\label{cl2}
Let $\lambda\ge \LS(\K)+\kappa(\K)$, and suppose $(<\lambda,\le
n+1)$-existence
and weak $(<\lambda,n)$-uniqueness hold, $n\ge 2$.
Let $\S_1\prec \S_2$ be independent stable $(\mu,n)$- and
$(\lambda,n)$-systems inside some
models $M^{n}_1$, $M^{n}_2$ respectively. Then there is $\bar M^{n}_2\succ
M^{n}_2$ and an
embedding $f:M^n_1 \to \bar M^n_2$
such that $f\restriction \Union _{M\in \S_1} M = \id$ and the system
$\S_1\cup \{f(M^n_1)\}\cup \S_2$
is a stable $(\lambda,n+1)$-system inside $\bar M^n_2$.
\end{lemma}

\begin{proof}
Iterate Lemma~\ref{cl1} $\lambda$-many times.
\end{proof}

\begin{proof}[Proof of Theorem~\ref{ex transfer}.]

Let $\S=\{M^{s}\mid s \in \Pm(n)\}\subset \K_\lambda$ be an (incomplete)
system of models. Our goal is to find a model $M^n$ and the coherent
embeddings
$f^s:M^s \to M^n$.

Take $\S_i:=\{M^s_i\mid i<\mu,s\in \Pm(n)\}$ a resolution of the system $\S$
such that for all $s\in \Pm(n)$ $\|M^s_i\|=\chi+|i|$ and $\S_i\prec \S_j$
for $i<j$.

For the base case, we just take a completion $M^n_0$ of the stable system
$\{M^s_0\mid s\in \Pm(n)\}$. Namely, we get a system of mappings
$f^s_0:M^s_0\to M^n_0$. It exists since we are assuming
$(\chi,n)$-existence.

Successor step. We have the model $M_i^n$, in which $f^s_i(M_i^s)$, $s\in
\Pm(n)$, form a stable
$n$-system. And from the resolution we have $M^s_{i+1}$ for $s\in \Pm(n)$,
$|s|=n-1$, where
$\{M^t_j \mid (\emptyset,i)\le (t,j)<(s,i+1)\}$, form a stable $n$-system
in size $\mu=\chi+|i|$.

By $(\mu,n+1)$-existence, we get $M_{i+1}^{n}$ and embeddings
$f_{i+1}^s:M_{i+1}^s\to M_{i+1}^n$ for $s\subset _{n-1}n$.
Now $(\mu,n+1)$-amalgamation also gives that $f^s_{i+1}\supset f^s_i$ for
$s\in \Pm(n)$.

For the limit step we simply take the union. Finally, the model
$M^{\lambda}_{n}$ is as needed.
\end{proof}

\begin{proof}[Proof of Theorem~\ref{uq transfer}.]
Let $\S^\ell$, $\ell=1,2$,
be stable $(\lambda,n)$-systems that are piecewise isomorphic.
We are constructing models $\bar N^\ell$ extending $M^\ell_n$ and a
$\K$-isomorphism
$\bar f:\bar N^1\to \bar N^2$ that extends all $f_s:M^1_s\to M^2_s$ for
$s\in \Pm(n)$.

By definability of forking we can find $M^\ell_{0,s}\prec M^\ell_s$
such that $\|M^\ell_{0,s}\|=\chi$ and $A^\ell_n\nonfork
^{N^\ell}_{A^\ell_{0,n}} M^\ell_{0,n}$.
Let $\{A^\ell_{i,n} \mid i<\lambda \}$ be an increasing continuous chain
whose union is $A^\ell_n$ and
$\|A^\ell_{i,n}\|=\chi+|i|$.

By induction on $i<\lambda$ build models $\bar N^\ell_i$, $\hat N^\ell_i$,
and isomorphisms $\bar f_i:\bar N^1_i \to \bar N^2_i$
such that
\begin{enumerate}
\item
$\bar N^\ell_{i}\succ N^\ell_{i}$;
\item
$\bar f_{i}\restriction A^1_{i,n}=\Union_{s\in \Pm(n)} f_s(M^1_{s,n})$;
\item
$\bar f_{i+1}$ extends $\bar f_i$;
\item
$A^\ell_n \nonfork^{\hat N^\ell_i}_{A^\ell_{i,n}} \bar N^\ell_i$ ($\|\hat
N^\ell_i\|=\lambda$).
\end{enumerate}

Begin with $i=0$. By weak $(\chi,n)$-uniqueness there is a model $\bar
N^2_0 \succ N^2_0$
and an embedding $f_0:N^1_0 \to \bar N^2_0$. By extension, we may assume
that
there is $\hat N^2_0\succ A^2_n, \bar N^2_0$ such that
$A^2_n \nonfork _{A^2_{0,n}}^{\hat N^2_0} \bar N^2_0$.
Let $\bar N^1_0\in \K$ be such that $\bar N^1_0\succ N^1_0$
and $\bar N^1_0$ is isomorphic to $\bar N^2_0$ via some $\bar f_0$ such that
$\bar f_0\restriction N^1_0 = f_0$.
Using extension again, we get $\hat N^1_0$ such that
$A^1_n \nonfork^{\hat N^1_0}_{A^1_{0,n}} \bar N^1_0$.

For $\alpha$ a limit ordinal, let $\bar N^\ell_\alpha:=\Union_{i<\alpha}
\bar N^\ell_i$;
$\hat N^\ell_\alpha:=\Union_{i<\alpha} \hat N^\ell_i$; and $\bar
f_\alpha:=\Union_{i<\alpha} \bar f_i$.
It is routine to check that (1)--(3) hold, and we need to establish (4). By
the induction hypothesis
and monotonicity, for all $i<\alpha$ we have
$A^\ell_n \nonfork^{\hat N^\ell_\alpha}_{A^\ell_{\alpha,n}} \bar N^\ell_i$.
So by continuity we get
$A^\ell_n \nonfork^{\hat N^\ell_\alpha}_{A^\ell_{\alpha,n}} \bar
N^\ell_\alpha$.

For the successor case, let $\mu:=\delta +|i|$. Let $N^\ell_{i+1}\succ
N^\ell_i$ be a $\K$-submodel
of $\hat N^\ell_i$ containing $A^\ell_{i+1,n}$; $\|N^\ell_{i+1}\|=\mu$.
By monotonicity, $A^\ell_{i+1,n}\nonfork^{N^\ell_{i+1}}_{A^\ell_{i,n}} \bar
N^\ell_i$,
so the system $\S^\ell_{i+1}:=A^\ell_{i+1,n}\cup A^\ell_{i,n}\cup \{\bar
N^\ell_i\}$ is a $(\mu,n+1)$-stable
system inside $N^\ell_{i+1}$. By weak $(\mu,n+1)$-uniqueness
($\S^\ell_{i+1}$, $\ell=1,2$ are
piecewise isomorphic), there is a model $\bar N^2_{i+1} \succ N^2_{i+1}$
and an embedding $f_{i+1}:N^1_{i+1} \to \bar N^2_{i+1}$ that extends the
``piecewise isomorphisms''
$f_{i+1,s}:M^1_{i+1,s}\cong M^2_{i+1,s}$ as well as the isomorphism $\bar
f_i$.
By $(<\lambda,n)$-existence and Lemma~\ref{cl2}, we may assume that
there is $\hat N^2_{i+1}\succ \hat N^2_{i}, \bar N^2_{i+1}$ such that
$A^2_n \nonfork^{\hat N^2_{i+1}}_{A^2_{i+1,n}} \bar N^2_{i+1}$.

Using Fact~\ref{pull-back}, we get $\bar N^1_{i+1}\in \K_\mu$ such that
$\bar N^1_{i+1}\succ N^1_{i+1}$
and $\bar N^1_{i+1}$ is isomorphic to $\bar N^2_{i+1}$ via some $\bar
f_{i+1}$ such that
$\bar f_{i+1}\restriction N^1_{i+1} = f_{i+1}$. Using
$(<\lambda,n)$-existence again, we get
$\hat N^1_{i+1}$ such that
$A^1_n \nonfork^{\hat N^1_{i+1}}_{A^1_{i+1,n}} \bar N^1_{i+1}$.
\end{proof}

%%% end section 3

\section{Three dimensional amalgamation}

A previous draft of this paper dealt with $n$-dimensional amalgamation
properties. In this section, we state \emph{a} definition of 3-dimensional
amalgamation, outline the proof of $(\lambda,\lambda^+)$-tameness from
$(\lambda,3)$-amalgamation, and finally show that
$(\lambda,3)$-amalgamation implies the weak $(\lambda,2)$-uniqueness
property.

The outline of the proof was presented by Rami
Grossberg in Bogot\'{a} model theory conference in the fall in
2003 and a preliminary version of this paper was posted on the web since December 2003.  In August 2005, weeks after we completed our proof we have learned that our idea of using 3-dimensional amalgamation was used to show directly that a certain natural class of structures is tame.
Using variants of 3-dimensional amalgamation  Villaveces-Zambrano in~\cite{ViZa} and Baldwin in~\cite{Ba} managed to
obtain tameness of certain abstract elementary classes arising naturally from Hrushovski's fusion of strongly minimal theories.  As the work of Villaveces-Zambrano and Baldwin is still in progress we suggest to the interested reader to consult them for their most recent results.

\begin{definition}
We say $\l \K,\nonfork\r$ has $(\lambda,3)$-amalgamation if for any system
of seven $\K_\lambda$-models $\{M_i\mid i<7\}$ such that
$M_1 \nonfork_{M_0}^{M_5} M_4$ and $M_2 \nonfork_{M_0}^{M_6} M_4$ and a
$\K$-embedding $f_0$, there
is a model $N^*$ and embeddings
$f:M_6\to N^*$ $g:M_5\to N^*$ and $h:M_3\to N^*$ such that the following diagram commutes.
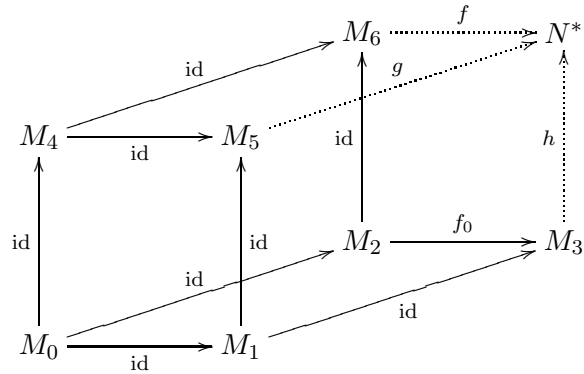
\[
      \xymatrix{\ar @{}
    &&& M_6\ar@{.>}[rr]^{f}  && N^*
\\
M_4\ar[rrru]^{\id}  \ar[rr]_{\id} &
    & M_5 \ar@{.>}[rrru]^{g}
\\ &&& M_2 \ar[uu]^{\id} \ar[rr]^{f_0} && M_3 \ar@{.>}[uu]^{h}
\\
M_0\ar[uu]^{\id} \ar[rrru]^{\id}  \ar[rr]_{\id} &
&  M_{1} \ar[uu]_{\id}  \ar[rrru]_{\id}
      }
\]
\end{definition}

\subsection{Tameness from 3-dimensional amalgamation}

Let $p, q\in \gaS(M)$ such that $M\in \K_\lambda$ and
$q\restriction N=p\restriction N$ for all $N\in\K_{<\lambda}$  enough to
show
that this condition implies $p=q$.

Suppose $(a_1,M,N^1)\in p$ and $(a_2,M,N^2)\in q$.  By the
$\mbox{LS}$-axiom pick
$\{N^\ell_\alpha\prec_{\K}N^\ell\mid \alpha<\lambda\}\subseteq
\K_{<\lambda}$
increasing and continuous resolutions of $N_\ell$ and
$\{M^\alpha\prec_{\K}M\mid \alpha<\lambda\}\subseteq \K_{<\lambda}$
increasing
resolution of $M$ such that
$M_\alpha\prec N^\ell_\alpha$, require that $a_\ell\in N^0_\ell$.

By the assumption there exist $N^*_0\succ_{\K}N^2_0$ of cardinality
less than $\lambda$, an amalgam of $N^1_0$ and $N^2_0$ over $M_0$
mapping $a_1$ to $a_2$.

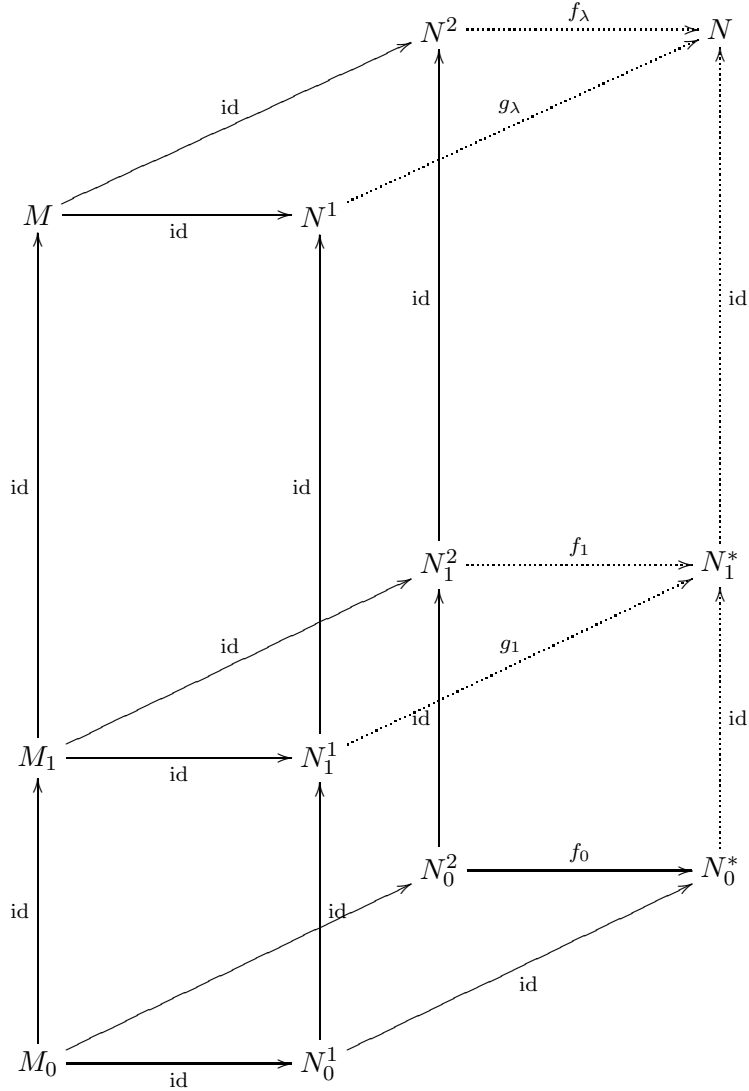
\[
     \xymatrix{\ar @{}   &&& & N^2  \ar@{.>}[rrr]^{f_\lambda}  &&& N\\  \\
     M
     \ar[rrr]_{\id} \ar[rrrruu]^{\id}  & & & N^1\ar@{.>}[rrrruu]^{g_\lambda}
\\
\\ \\
\\ &&&& N^2_1 \ar[uuuuuu]^{\id} \ar@{.>}[rrr]^{f_1} &&&
N^*_1\ar@{.>}[uuuuuu]_{\id}\\ \\
M_1\ar[uuuuuu]^{\id} \ar[rrrruu]^{\id}  \ar[rrr]_{\id} &
& & N^1_1\ar[uuuuuu]^{\id}  \ar@{.>}[rrrruu]^{g_1}
\\ &&&& N^{2}_0 \ar[uuu]^{\id} \ar[rrr]^{f_0} &&& N^*_0
\ar@{.>}[uuu]_{\id}\\ \\
M_0\ar[uuu]^{\id} \ar[rrrruu]^{}  \ar[rrr]_{\id} &
& & N_0^{1} \ar[uuu]_{\id}  \ar[rrrruu]_{\id}
     }
\]

Clearly it is enough to find $N_1^*\in \K_{<\lambda}$ a
$\prec_{\K}$-extension
of $N^*_0$ and $\K$-embeddings
$g_1:N_1^1 \rightarrow N^*_1$ and $f_1:N^2_1\rightarrow N^*_1$ such that
the above diagram commutes.  Continuing this by induction on
$\alpha<\lambda$
gives $N^*_\alpha$ and $f_\alpha,g_\alpha$ such that the above diagram
commutes.

Let $N:=\bigcup_{\alpha<\lambda}N_\alpha^*$ and
$f_\lambda:=\bigcup_{\alpha<\lambda}f_\alpha$.
Since $f_\lambda(a_1)=a_2$ we have that $p=q$.

\subsection{3-amalgamation implies weak 2-uniqueness}
Rather than formalize the above argument, we show that $(\lambda,3)$-amalgamation
is a strong enough assumption, so that a particular case of
$(\lambda,3)$-amalgamation implies weak $(\lambda,2)$-uniqueness.

\begin{proposition}
Suppose $\K$ has $(\lambda,3)$-amalgamation, where all the given embeddings are
identity (in other words, we allow only systems where $f_0=\id$). Then weak
$(\lambda,2)$-uniqueness holds.
\end{proposition}

\begin{proof}
Take $M_0\prec M_1,M_2 \in \K_\lambda$; let $N,N'\in \K_\lambda$ both contain
$M_1\cup M_2$ and suppose that $M_1\nonfork_{M_0}^{N,N'} M_2$.
By Proposition~\ref{wu=failure of non uniq}, it is enough to show that
$N$, $N'$ can be amalgamated over $M_1\cup M_2$.

We have the following diagram:
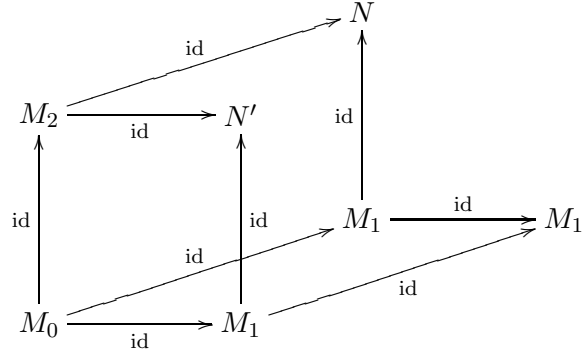
\[
      \xymatrix{\ar @{}
    &&& N  &&
\\
M_2\ar[rrru]^{\id}  \ar[rr]_{\id} &
    & N'
\\ &&& M_1 \ar[uu]^{\id} \ar[rr]^{\id} && M_1
\\
M_0\ar[uu]^{\id} \ar[rrru]^{\id}  \ar[rr]_{\id} &
&  M_{1} \ar[uu]_{\id}  \ar[rrru]_{\id}
      }
\]
By $(\lambda,3)$-amalgamation, there is a model $N^*$ and embeddings
$f:N\to N^*$
$g:N'\to N^*$ and $h:M_1\to N^*$ such that the following diagram commutes.
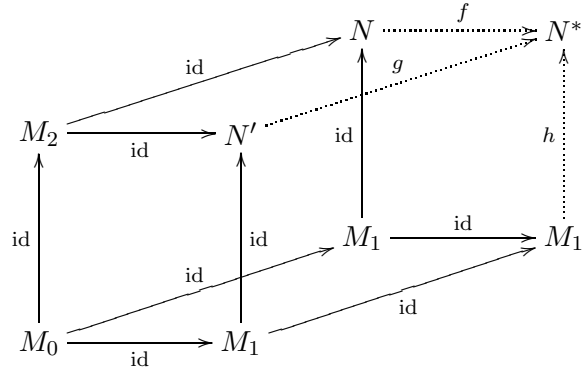
\[
      \xymatrix{\ar @{}
    &&& N\ar@{.>}[rr]^{f}  && N^*
\\
M_2\ar[rrru]^{\id}  \ar[rr]_{\id} &
    & N' \ar@{.>}[rrru]^{g}
\\ &&& M_1 \ar[uu]^{\id} \ar[rr]^{\id} && M_1 \ar@{.>}[uu]^{h}
\\
M_0\ar[uu]^{\id} \ar[rrru]^{\id}  \ar[rr]_{\id} &
&  M_{1} \ar[uu]_{\id}  \ar[rrru]_{\id}
      }
\]
Now $f(M_2)=g(M_2)$, and $f(M_1)=h(M_1)=g(M_1)$, so $f(M_i)=g(M_i)$ for
$i=0,1,2$.
Thus $N^*$ is an amalgam of $N$ and $N'$ over $M_1\cup M_2$.
\end{proof}

%%%%%%%%%%%%%%%%%%%%%%%%%%%%%%%%%%%%%%%%%%%%%%%%%%%%%%%%%%%%%%%%%
%%%%%%%%%%%%%%%%%%%%%%

\end{document}